\documentclass[a4paper,12pt]{amsart}
\usepackage{color}
\usepackage{amsmath,amssymb,amscd,amsfonts,amsthm,stmaryrd}
\usepackage[mathscr]{eucal}
\usepackage{xcolor}
\usepackage{graphicx}
\usepackage{tikz-cd}
\usepackage{adjustbox}

\usepackage{nicefrac}

\usepackage{color}

\numberwithin{equation}{section}

\def\length{\operatorname{length}}

\def\diam{\operatorname{diam}}

\def\pol{\operatorname{pol}}

\def\acts{\curvearrowright}
\def\D{\partial}

\def\R{\mathbb R}
\def\Z{\mathbb Z}

\def\N{\mathbb N}
\def\Ha{\mathcal H}
\def\HH{\mathbb{H}}

\newcommand{\F}{F}
\newcommand{\G}{\Theta}
\newcommand{\Fi}{\F_\infty}
\newcommand{\Gi}{\G_\infty}

\def\al{\alpha}

\def\ka{\kappa}

\def\eps{\epsilon}
\def\ga{\gamma}
\def\Ga{\Gamma}

\def\la{\lambda}
\def\La{\Lambda}

\def\om{\omega}

\def\si{\sigma}
\def\Si{\Sigma}

\def\tits{\partial_{T}}
\def\geo{\partial_{\infty}}
\def\etits{\partial^{ess}_{T}}

\newcommand{\bI}{{\mathbf I}}
\newcommand{\bZ}{{\mathbf Z}}
\newcommand{\M}{{\mathbf M}}

\newcommand{\cD}{{\mathscr D}}
\newcommand{\cF}{{\mathscr F}}

\newcommand{\cP}{{\mathscr P}}

\newcommand{\on}{\:\mbox{\rule{0.1ex}{1.2ex}\rule{1.1ex}{0.1ex}}\:}

\def\acts{\curvearrowright}

\def\spt{\operatorname{spt}}
\def\Ant{\operatorname{Ant}}

\def\<{\langle}
\def\>{\rangle}

\newcommand{\bb}[1]{\llbracket #1\rrbracket} 
\newcommand{\loc}{\text{\rm loc}}

\theoremstyle{plain}
\newtheorem{thm}{Theorem}[section]
\newtheorem{lem}[thm]{Lemma}
\newtheorem{prop}[thm]{Proposition}
\newtheorem{cor}[thm]{Corollary}
\newtheorem{slem}[thm]{Sublemma}

\newtheorem{introthm}{Theorem}

\newtheorem{introcor}[introthm]{Corollary}

\newtheorem{conj}{Conjecture}

\newtheoremstyle{named}{}{}{\itshape}{}{\bfseries}{.}{.5em}{\thmnote{#3} #1}
\theoremstyle{named}
\newtheorem*{namedlemma}{Lemma}

\theoremstyle{definition}
\newtheorem{dfn}[thm]{Definition}

\theoremstyle{remark}
\newtheorem{rem}[thm]{Remark}

\newcommand{\bcl}{\begin{claim}}
\newcommand{\ecl}{\end{claim}}
\newcommand{\bcor}{\begin{cor}}
\newcommand{\ecor}{\end{cor}}
\newcommand{\bdfn}{\begin{dfn}}
\newcommand{\edfn}{\end{dfn}}
\newcommand{\ben}{\begin{enumerate}}
\newcommand{\bit}{\begin{itemize}}
\newcommand{\blem}{\begin{lem}}
\newcommand{\bslem}{\begin{slem}}
\newcommand{\bprop}{\begin{prop}}
\newcommand{\bthm}{\begin{thm}}
\newcommand{\een}{\end{enumerate}}
\newcommand{\eit}{\end{itemize}}
\newcommand{\elem}{\end{lem}}
\newcommand{\eslem}{\end{slem}}
\newcommand{\eprop}{\end{prop}}
\newcommand{\ethm}{\end{thm}}

\usepackage{hyperref}

\begin{document}

\title[Rank Rigidity]{Rank Rigidity for CAT(0) spaces without 3-flats}
\author{Stephan Stadler}

\newcommand{\Addresses}{{\bigskip\footnotesize
\noindent Stephan Stadler,
\par\nopagebreak\noindent\textsc{Max Planck Institute for Mathematics, Vivatsgasse 7, 53111 Bonn, Germany}
\par\nopagebreak
\noindent\textit{Email}: \texttt{stadler@mpim-bonn.mpg.de}

}}

%\subjclass[2010]{53C20, 53C23, 58E20}

%\keywords{Non-positive curvature, conformal change, minimal disc, harmonic maps}

\begin{abstract}
If a group $\Ga$ acts geometrically on a CAT(0) space $X$
without 3-flats, then either $X$ contains a $\Ga$-periodic geodesic which does not bound a flat half-plane,
or else $X$ is a rank 2 Riemannian symmetric space, a 2-dimensional Euclidean building or non-trivially splits 
as a metric product. Consequently all such groups satisfy a strong form of the Tits Alternative.
\end{abstract}

\maketitle

\tableofcontents

%TO DO:
%\begin{itemize}
%\item abstract
%\item organization
%\item rewrite every statement according to rank
	%\item acknowledgement
	%\item pictures
	%\item intro sentence to sections
	%\item check for defining signs
%\item unify setting within section
%\item appendix
%
%\end{itemize}

\section{Introduction}
\subsection{Main results}
In \cite{G_asym} Gromov proposed a rough classification of periodic CAT(0) spaces modulo hyperbolic ones.
It vaguely states that these spaces decompose into pieces such that every piece either displays
almost hyperbolic behavior or else has uniformly distributed flats. One may want to think of
a non-positively curved compact 3-manifold and its decomposition into Seifert and atoroidal pieces.
Each Seifert piece admits an $\R^3$- or 
$\HH^2\times\R$-structure whereas the atoroidal pieces admit 
$\HH^3$-structures.
Gromov's expectation roots in the general principle that the failure of hyperbolic behavior is caused by (unbounded parts of) flats --
isometric embeddings of higher dimensional Euclidean spaces. For instance, an axial isometry exhibits north-south dynamics
at infinity unless one (and then every) of its axes bounds a flat half-plane.

One concrete interpretation of Gromov's vision is provided by the following two conjectures due to Ballmann and Buyalo \cite{BaBu_periodic}
which classify CAT(0) spaces by means of their {\em rank}.
A complete geodesic is said to have {\em rank 1} if it does not bound a flat half-plane.
Accordingly, a CAT(0) space is called {\em rank 1} if it contains a rank 1 geodesic,
otherwise we say it has {\em higher rank}.

Recall that if a group $\Ga$ acts on a CAT(0) space $X$, then its {\em limit set} $\La(\Ga)$ is the set of accumulation points  
of a $\Ga$-orbit in the ideal boundary $\geo X$.

\begin{conj}[Closing Lemma]\label{conj_cl}
Let $X$ be a locally compact CAT(0) space which admits a  properly discontinuous action $\Ga\acts X$ with $\La(\Ga)=\geo X$.
If $X$ contains a geodesic of rank 1, then it also contains a $\Ga$-periodic
geodesic of rank 1. 
\end{conj}
 
%Conjecture~\ref{conj_cl} holds under the assumption that $\Ga$ satisfies the duality condition of Chen 
%and Eberlein, compare \cite[Chapter~III]{ballmannbook}. In particular, Conjecture~\ref{conj_cl} holds in the case where 
%$X$ is a Hadamard manifold and $\Ga$ is a properly discontinuous group of isometries of $X$ 
%with finite covolume. Conjecture~\ref{conj_cl} also holds under the assumption that $X$ is a piecewise 
%smooth complex of dimension two and $\Ga$ acts {\em geometrically} -- properly discontinuously with compact quotient, see \cite{BaBr_orbi}.
%Moreover, Conjecture~\ref{conj_cl} holds for finite-dimensional CAT(0) cube complexes \cite{CS_rr}.

Conjecture~\ref{conj_cl} is known to hold in each of the following settings.

\begin{itemize}
	\item If $\Ga$ satisfies the duality condition of Chen 
and Eberlein, compare \cite[Chapter~III]{ballmannbook}.
\item If $X$ is a piecewise 
smooth complex of dimension two and $\Ga$ acts {\em geometrically} -- properly discontinuously with compact quotient, see \cite{BaBr_orbi}.
\item If $X$ is a finite-dimensional CAT(0) cube complex and $\Ga$ acts geometrically \cite{CS_rr}.
\end{itemize}

\begin{conj}[Diameter Rigidity]\label{conj_rr}
Let $X$ be a locally compact and geodesically complete CAT(0) space which admits a  properly discontinuous action $\Ga\acts X$ with $\La(\Ga)=\geo X$.
If $X$ has higher rank, then $X$ is a Riemannian symmetric space, a Euclidean building or non-trivially splits as a
metric product.
\end{conj}

Recall that in a CAT(0) space $X$, every complete geodesic bounds a flat half-plane if and only if its Tits boundary $\tits X$
has diameter equal to $\pi$, compare Section~\ref{sec_catzero}. This explains the name of the conjecture.

The following special cases of Conjecture~\ref{conj_rr} are known to hold. 

\begin{itemize}
	\item If $X$ is a Hadamard manifold and $\Ga$ satisfies the duality condition \cite{B_higher,BS_higher,EH_diff}; for a comprehensive proof, see Chapter IV in 
\cite{ballmannbook}.
\item If $X$ is a homogeneous Hadamard manifold \cite{Heber}. 
\item If $X$ is a piecewise 
smooth complex of dimension two or a piecewise Euclidean complex of dimension 3 and 
$\Ga$ acts geometrically \cite{BaBr_orbi,BB_rr}.
\item If $X$ is a finite-dimensional
CAT(0) cube complex and $\Ga$ acts geometrically \cite{CS_rr}.
\end{itemize}

We confirm Conjectures~\ref{conj_cl} and~\ref{conj_rr} for CAT(0) spaces with geometric group actions and without 3-flats.
By \cite[Theorem~C]{Kleiner}, the absence of 3-flats is equivalent to $X$ not containing 3-dimensional quasi-flats.

\begin{introthm}\label{thm_mainA}
Let $\Ga$ be a group acting geometrically on a locally compact CAT(0) space without $3$-flats.
If $X$ contains a  geodesic of rank 1, then it also contains a
 $\Ga$-periodic geodesic of rank 1.
\end{introthm}

\begin{introthm}\label{thm_mainB}
Let $\Ga$ be a group acting geometrically on a locally compact and geodesically complete CAT(0) space without $3$-flats.
If $X$ has higher rank,
then $X$ is a Riemannian symmetric space, a Euclidean building or non-trivially splits as a metric product.
\end{introthm}

A consequence of Theorems~\ref{thm_mainA} and~\ref{thm_mainB} is that all such groups satisfy the Tits Alternative:

\begin{introcor}\label{cor_tits}
Let $\Ga$ be a group acting geometrically on a locally compact and geodesically complete CAT(0) space without $3$-flats.
Then either $\Ga$ contains a free non-abelian subgroup or else $X$ is flat and $\Ga$ is a Bieberbach group.
\end{introcor}

Another application of Theorem~\ref{thm_mainA} gives  the following.

\begin{introcor}\label{cor_main}
Let $X$ and $X'$ be irreducible locally compact geodesically complete CAT(0) spaces. Suppose that $X$ does not contain a 3-flat.
If the same group $\Ga$ acts geometrically on $X$ and $X'$,
then either $X$ and $X'$ are isometric after possibly rescaling or there exists an element $\ga\in\Ga$ with $\ga$-axes $c\subset X$
and $c'\subset X'$ such that neither of them  bounds a flat half-plane.
\end{introcor}

%We also obtain
%
%\begin{introcor}\label{cor_qi}
%Let $X$ be an irreducible locally compact geodesically complete CAT(0) space with a geometric group action. 
%If $X$ is quasi-isometric to a rank 2 Riemannian symmetric space of Euclidean building $X_{model}$, then
%$X$ is isometric to a rescaling of $X_{model}$.
%\end{introcor}
%\medskip
%
%Note that this can also be deduced from \cite{KL_induced} in combination with \cite{Leeb}.

The present article is part of a series motivated by the Diameter Rigidity Conjecture (sometimes also referred to
as Higher Rank Rigidity). The main results obtained in the other parts are as follows.

\bthm[{\cite[Theorem~A]{St_rrI}}]\label{thm_rrI}
Let $X$ be a locally compact  CAT(0) space whose Tits boundary has dimension $n-1\geq 1$. 
Suppose that every geodesic in $X$ lies in an $n$-flat. 
If $X$ contains a periodic $n$-flat, then $X$ is a Riemannian symmetric space or a Euclidean building, or  $X$ non-trivially splits as a metric product.
\ethm

\bthm[{\cite[Main Theorem]{St_rrII}}]\label{thm_rrII}
Let $X$ be a locally compact CAT(0) space  with a geometric group action $\Ga\acts X$.  
Suppose that there exists $n\geq 2$ such that every geodesic in $X$ lies in an $n$-flat. 
If $X$ contains a periodic $n$-flat which does not bound a flat half-space, then $X$ is a Riemannian symmetric space or a Euclidean building or  $X$
non-trivially splits as a metric product.
\ethm

\subsection{Organization}

In Section~\ref{sec_pre}, we introduce notation and recall background on the geometry of spaces with upper curvature bounds.
Apart from the standard material, we introduce the notion {\em essential Tits boundary} of a CAT(0) space without 3-flats.
In Section~\ref{sec_cl} we proof the Closing Lemma by showing that if a group $\Ga$ acts geometrically on a CAT(0) space without 3-flats
and preserves a proper closed subset at infinity, then the space has higher rank (Theorem~\ref{thm_diambound}).
Section~\ref{sec_rank_2} builds up to a proof of Theorem~\ref{thm_mainB} assuming our key technical result, the \hyperref[lem_key_tech]{Half-Plane Lemma}.
We proceed by first showing that the essential Tits boundary is a  spherical building and then argue that the essential Tits boundary
already agrees with the ordinary Tits boundary.

Appenix~\ref{sec_app} is devoted to geometric measure theory with the aim to provide a proof of the \hyperref[lem_key_tech]{Half-Plane Lemma}.
We also prove some folklore results, including monotonicity and volume rigidity of minimizing currents in CAT(0) spaces.
However, our approach is to find a streamlined path to  the \hyperref[lem_key_tech]{Half-Plane Lemma} and we do not try to prove
required ingredients in their most general form.

\subsection{Acknowledgments}
It's my pleasure to thank several people for their support. 
I want to thank Alexander Lytchak for  answering many questions, for reading a late version of this article and 
for his critical and very helpful feedback. 
I want to thank Bruce Kleiner for several inspiring discussions.
This paper draws heavily from his ideas and his work.
I want to thank Bernhard Leeb for valuable comments.
 %I would like to use this opportunity to thank Bernhard Leeb and Alexander Lytchak for 
%pointing out several mistakes in a late version of this article.
%Their insightful comments led to substantial improvements. 
I also want to thank Anton Petrunin for a helpful discussion on
Proposition~\ref{prop_mon}. 
I was supported by DFG grant SPP 2026.

%\subsection{Acknowledgments}
%It's my pleasure to thank Alexander Lytchak for carefully reading a late version of this paper,
%pointing out several mistakes,
%and making insightful comments that led to substantial improvements. 
%I also want to thank Anton Petrunin for a helpful discussion on
%Proposition~\ref{prop_mon}. 
%I was supported by DFG grant SPP 2026.

\section{Preliminaries}\label{sec_pre}

General references for this section are \cite{AKP, Ballmann, BH, KleinerLeeb}.

\subsection{Metric spaces}
Euclidean $n$-space with its flat metric will be denoted by $\R^n$. The unit sphere $S^{n-1}\subset\R^n$ equipped with the induced metric will be referred to as a
{\em round sphere}. Its intersection with a half-space $\R^{n-1}\times[0,\infty)$ is called a {\em round hemisphere}.
We denote  the distance between two points $x$ and $y$ in a metric space $X$ by $|x,y|$.
If $A\subset X$ denotes a subset, then $|x,A|$ refers to the greatest lower bound for distances from points in $A$ to $x$.
%For subsets $A, A'\subset X$ we denote the Hausdorff (pseudo-) distance by $|A,A'|_H$.  
For $x\in X$ and $r>0$, we denote by $B_r(x)$ and $\bar B_r(x)$ the open and closed $r$-ball around $x$, respectively.
Similarly, $N_r(A)$ and $\bar N_r(A)$ denote the open and closed $r$-neighborhood of a subset $A\subset X$, respectively.
Moreover, $S_r(x)$ denotes the $r$-sphere around $x$ and by $\dot B_r(x)$ we denote the punctured $r$-ball $B_r(x)\setminus\{x\}$.
A \emph{geodesic}
is an isometric embedding of an interval. It is called a {\em geodesic segment}, if it is compact.
The {\em endpoints} of a geodesic segment $c$ are denoted by $\D c$.
A geodesic segment $c$ {\em branches} at an endpoint $y\in \D c=\{x,y\}$, if there are geodesics $c^\pm$ starting in $x$ which  strictly contain $c$
and such that $c^-\cap c^+=c$. The point $y$ is then called a {\em branch point}.

 A \emph{triangle} is a union of three geodesics connecting three points.
$X$  is \emph{a  geodesic metric space} if
any pair of   points of $X$
is connected by a geodesic.
It is \emph{geodesically complete} if every geodesic segment is contained in a complete local geodesic.

\subsection{Spaces with an upper curvature bound}

For $\kappa \in \R $, let $D_{\kappa}  \in (0,\infty] $ be the diameter of the  complete, simply connected  surface $M^2_{\kappa}$
of constant curvature $\kappa$. A complete  metric space is called a \emph{CAT($\kappa$) space}
if any pair of its points with distance less than $D_{\kappa}$ is connected by a geodesic and if
 all triangles
with perimeter less than $2D_{\kappa}$
are not thicker than
 the \emph{comparison triangle} in $M_{\kappa} ^2$. In particular, geodesics between points of distance less than $D_{\kappa}$
are unique. Hence, $X$ is a CAT($\ka$) space, then we can define for every $p\in X$ and every subset $A\subset B_{D_{\ka}}(p)$ 
the {\em geodesic cone $C_p(A)$}.

\subsection{Directions and angles}
For any CAT($\kappa$) space  $X$,
the angle between each pair of geodesics starting at the same point
is well defined. 
 Let $x,y,z$ be three points at pairwise distance less than $D_{\kappa}$ in a CAT($\kappa$) space $X$.
Whenever $x\neq y$, the geodesic between $x$ and $y$ is unique and will be denoted
by $xy$.   For $y,z \neq x$, the angle at $x$ between $xy$ and $xz$
will be denoted by $\angle_x(y,z)$. It is defined by 
\[\angle_x(y,z)=\lim\limits_{y',z'\to x}\tilde\angle_x(y',z')\]
where $\tilde\angle_x(y',x')$ denotes the angle of the comparison triangle at the vertex corresponding to $x$;
and the points $y',z'$ converge to $x$ along the geodesic $xy$ and $xz$, respectively.
A comparison argument shows that this is well defined as a monotonic limit.
The \emph{space of directions} or {\em link}
 at a point  $x\in X$ is the completion of the space of equivalence classes of geodesic germs at $x$ with respect to the angle metric.  
The resulting space $(\Si_x X,\angle)$
is a CAT(1) space. Its elements are called {\em directions (at $x$)}. 

\subsection{Dimension}
 A natural notion of dimension $\dim (X)$ for a CAT($\ka$) space $X$ was
  introduced by
Kleiner in \cite{Kleiner}, originally referred to as {\em geometric dimension}.
It vanishes precisely when the space is discrete.
In general, it is defined inductively:
\[\dim (X)= \sup _{x\in X} \{ \dim (\Sigma _xX) +1 \}.\]
For instance, in a 1-dimensional CAT($\ka$) space every link is discrete. 
The geometric dimension coincides with the supremum of 
topological dimensions\footnote{Here topological dimension corresponds to Lebesgue covering dimension.} of compact subsets in $X$ \cite{Kleiner}. 
If the dimension of $X$ is finite, then $\dim(X)$ agrees with the largest number $k$
such that $X$ admits a bilipschitz embedding of an open subset in $\R^k$ \cite{Kleiner}.
The dimension of a locally compact and geodesically complete space is finite and agrees with the topological dimension 
as well as the Hausdorff dimension \cite{OT_cba, LN_gcba}.

\subsection{CAT(1) spaces}

For two CAT(1) spaces $Z_1$ and $Z_2$ we denote by $Z_1\circ Z_2$ their {\em spherical join}.
It is a CAT(1) space of diameter $\pi$.
Note that every round subsphere $\si'$ in a round sphere $\si$ yields a join decomposition $\si=\si'\circ\si''$
where $\dim(\si)=\dim(\si')+\dim(\si'')+1$. Also, every round hemisphere $\tau$ decomposes as $\tau=\D\tau\circ\zeta$
where $\zeta$ denotes the center of $\tau$. A CAT(1) space $Z$ is called {\em irreducible}, if it does not admit a non-trivial
spherical join decomposition.

 %A subset $C$ in a CAT(1) space $Z$ is called {\em $\pi$-convex}
%if for any pair of points $x, y\in Z$ at distance less than $\pi$ the unique geodesic $xy$ is contained in $C$.
%If $C$ is closed, then it is CAT(1) with respect to the induced metric.
%For instance, a ball of radius at most $\frac{\pi}{2}$ is $\pi$-convex.
%Let $C\subset Z$ be a closed convex subset with radius $\rad(C)\geq \pi$. Then we define the set of {\em poles} for $C$ by
%\[\pol(C):=\{\eta\in Z|\ d(\eta,\cdot)|_{C}\equiv\frac{\pi}{2}\}.\]
%If $\diam(C)>\pi$, then $C$ has no pole. If $\diam(C)=\rad(C)=\pi$, then $\pol(C)$ is closed and convex.
%The convex hull of $C$ and $\pol(C)$ is canonically isometric to $C\circ\pol(C)$.
Recall that two points in a CAT(1) space $Z$ are called {\em antipodes}, if their distance is at least $\pi$.
In case $C$ consists of a pair of {\em antipodes} $\xi^\pm$, $|\xi^-,\xi^+|=\pi$, then the convex hull of
$\{\xi^-,\xi^+\}$ and $\pol(\{\xi^-,\xi^+\})$ is isometric to a spherical suspension of $\pol(\{\xi^-,\xi^+\})$.
For a subset $M\subset Z$ we denote by $\Ant(M)\subset Z$ the set of antipodes of $M$.

\bdfn\label{def_reg_point}
Let $Z$ be a CAT(1) space of dimension $k$. We call a point $\xi\in Z$ {\em regular}, if
it has a neighborhood homeomorphic to an open set in $\R^k$. 
\edfn

\blem[{\cite[Lemma~2.1]{BL_building}}]\label{lem_antintop}
	Let $Z$ be a CAT(1) space of dimension $n$. If $\si\subset Z$ is a round $n$-sphere, then every point $\xi\in Z$ has an antipode in $\si$.
\elem

%\proof
%We prove the claim by induction on $n$. If $n=0$, then $Z$ is discrete and there is nothing to show.
%Now suppose the claim holds for all CAT(1) spaces of dimension at most $n-1$.
%Let $\eta$ be a  point in an $n$-dimensional CAT(1) space $Z$ and let $\si\subset Z$ be a round $n$-sphere.
%Choose a point $\zeta$ in $\si$. If it is not an antipode of $\eta$, then denote by $v\in\Si_\zeta Z$
%the direction towards $\eta$. Since $\dim(\Si_\zeta Z)\leq\dim(Z)-1$, we find an antipode $\hat v\in\Si_\zeta\si$ of $v$ by induction.
%Now we can extend the geodesic $\eta\zeta$ inside $\si$ in the direction $\hat v$ up to an antipode $\hat\eta$ of $\eta$.
%\qed
%\medskip

This implies the following well-known criterion for geodesic completeness.

\blem[{\cite[Lemma~2.4]{St_rrI}}]\label{lem_gc}
	Let $Z$ be a CAT(1) space of dimension $n$. If every point $\xi\in Z$ is contained in a round $n$-sphere in $Z$,
	then $Z$ is geodesically complete.
\elem
%
%\proof
%If $\xi$ is a point in $Z$ contained in a round $n$-sphere $\si\subset Z$, 
%then any local geodesic ending in the point $\xi$ can be extended beyond  $\xi$
%by a geodesic in $\si$, since any direction in $\Si_\xi Z$ has an antipode in $\Si_\xi\si\cong S^{n-1}$ by Lemma~\ref{lem_antintop}.
%\qed

\subsection{CAT(0) spaces}\label{sec_catzero}

The {\em ideal boundary} of a CAT(0) space $X$, equipped with the cone topology, is denoted by $\geo X$. 
If $X$ is locally compact, then $\geo X$ is compact. If $X$ is a Hadamard $n$-manifold, then $\geo X$ is homeomorphic to an $(n-1)$-sphere.  
The {\em Tits boundary} of $X$ is denoted by $\tits X$, it is the ideal boundary equipped with 
the Tits metric $|\cdot,\cdot|_T$. Recall that the Tits metric is the intrinsic metric associated to the {\em Tits angle}. For a point $p$ in $X$
and ideal points $\xi,\eta$ in $\geo X$ the Tits angle is defined by
\[\angle_T(\xi,\eta):=\lim\limits_{x\to\xi, y\to\eta}\tilde\angle_p(x,y).\]
 The Tits boundary of a CAT(0) space is a CAT(1) space. It decodes the intersection patterns of asymptotically flat subspaces.
For instance, the Euclidean $n$-space has the round $(n-1)$-sphere as Tits boundary, while any hyperbolic space has a discrete Tits boundary.

A subset in a CAT(0) space is convex, if it contains the geodesic between any pair of its points.

If $C$ is a closed convex subset, then it is CAT(0) with respect  to the induced metric. In this case, $C$ admits a 1-Lipschitz retraction $\pi_C:X\to C$.
If $C_1$ and $C_2$ are closed convex subsets, then the distance function $d(\cdot,C_1)|_{C_2}$ is convex; and constant if
and only if $\pi_{C_1}$ restricts to an isometric embedding on $C_2$. We call $C_1$ and $C_2$ {\em parallel}, $C_1\| C_2$, if
and only if  $d(\cdot,C_1)|_{C_2}$ and $d(\cdot,C_2)|_{C_1}$ are constant.
Let $Y\subset X$ be a geodesically complete closed convex subset. Then we define the {\em parallel set} $P(Y)$
as the union of all closed convex subsets parallel to $Y$. The parallel set is closed, convex and splits canonically as a metric product
\[P(Y)\cong Y\times CS(Y)\]
where the {\em cross section} $CS(Y)$ is a closed convex subset.
The Tits boundary is given by
\[\tits P(Y)\cong\tits Y\circ\pol(\tits Y).\]
If $X_1$ and $X_2$ are CAT(0) spaces, then their metric product $X_1\times X_2$ is again a CAT(0) space.
We have $\tits (X_1\times X_2)=\tits X_1\circ \tits X_2$ and $\Si_{(x_1,x_2)}(X_1\times X_2)=\Si_{x_1} X_1\circ \Si_{x_2} X_2$.  
If $X$ is a geodesically complete CAT(0) space, then any join decomposition of $\tits X$ is induced by a metric product decomposition of $X$ \cite[Proposition~2.3.7]{KleinerLeeb}.
A CAT(0) space $X$ is called {\em irreducible}, if it does not admit a non-trivial splitting  as a metric product.

The Tits boundary of a closed convex subset $Y\subset X$ embeds canonically $\tits Y\subset\tits X$.
If two closed convex subsets $Y_1$ and $Y_2$ intersect in $X$, then 
\[\tits(Y_1\cap Y_2)=\tits Y_1\cap \tits Y_2.\]

Let $X_1$ and $X_2$ be CAT(0) spaces containing closed convex subsets $C_1$ and $C_2$, respectively.
If there is an isometry $f:C_1\to C_2$, then by \cite[8.9.1 Reshetnyak's gluing theorem]{AKP}, the glued space
\[X_1\cup_f X_2:=(X_1\dot\cup X_2)/x_1\sim f(x_1)\]
is CAT(0) with respect to the induced length metric.
In particular, for every CAT(0) space $X$ which contains a closed convex subset $C$ we can construct a CAT(0) space $\hat X$, the {\em double}
of $X$ along $C$, by
\[\hat X=X^-\cup_C X^+\]
where $X^\pm$ are two isometric copies of $X$. The Tits boundary of the double is given by
\[\tits\hat X=\tits X^-\cup_{\tits C} \tits X^+.\]

For points $\xi\in\geo X$ and $x\in X$ we denote the {\em Busemann function centered at $\xi$ based at $x$} by $b_{\xi,x}$.
If $\rho:[0,\infty)\to X$ denotes the geodesic ray asymptotic to $\xi$ and with $\rho(0)=x$, then
\[b_{\xi,x}(y)=\lim\limits_{t\to\infty}(|y,\rho(t)|-t).\]
It is a 1-Lipschitz convex function whose negative gradient at a point $y\in X$ is given by $\log_y(\xi)$. 
We denote the {\em horoball} centered at a point $\xi\in\geo X$ and based at the point $x\in X$ by 
\[HB(\xi,x):=b^{-1}_{\xi,x}((-\infty,0]).\]
It is a closed convex subset with 
\[\tits HB(\xi,x)=\bar B_{\frac{\pi}{2}}(\xi)\subset\tits X.\]

A {\em $n$-flat} $F$ in a CAT(0) space $X$ is a closed convex subset isometric to $\R^n$.
In particular, $\tits F\subset\tits X$ is a round $(n-1)$-sphere.
On the other hand, if $X$ is locally compact and $\si\subset \tits X$ is a round  $(n-1)$-sphere, then either there exists an
$n$-flat $F\subset X$ with $\tits F=\si$, or there exists a round $n$-hemisphere $\tau^+\subset\tits X$
with $\si=\D\tau^+$ \cite[Proposition~2.1]{Leeb}. Consequently, if $\tits X$ is $(n-1)$-dimensional, then any round $(n-1)$-sphere in $\tits X$
is the Tits boundary of some $n$-flat in $X$. Moreover, a round $n$-hemisphere $\tau^+\subset\tits X$ bounds a flat $(n+1)$-half-space in $X$ if and only if
its boundary $\D\tau^+$ bounds an $n$-flat in $X$. A {\em flat ($n$-dimensional) half-space} $H\subset X$ is a closed convex subset isometric to a Euclidean half-space $\R_+^n$.
Its boundary $\D H\subset X$ is an $(n-1)$-flat and its Tits boundary is a round $(n-1)$-hemisphere $\tits H\subset\tits X$.  Flat half-spaces will play a certain role in our arguments later, and we agree to denote them by $H$
and their boundaries by $\D H=h$.

We define the {\em parallel set} of a round sphere $\si\subset \tits X$ as 
\[P(\si)=P(F)\]
where $F$ is a flat in $X$ with $\tits F=\si$, if such a flat exists.

\subsection{The essential Tits boundary}

\bdfn
For a CAT(0) space without 3-flats, we define its {\em essential Tits boundary} $\etits X$ as the subset of $\tits X$ given by the union of all simple closed local geodesics.
\edfn

Note that if $X$ is a CAT(0) space with isolated flats and 1-dimensional Tits boundary which supports a geometric group action, then $\etits X$ is dense in $\geo X$ but not closed.
However, we do have the following.

\blem\label{lem_diameter_etits}
Let $X$ be a CAT(0) space with $1$-dimensional Tits boundary which contains a 2-flat. If the diameter of $\tits X$ is equal to $\pi$, then 
$\etits X\subset\tits X$ is a convex subset of diameter $\pi$. In particular, $\etits X$ is a geodesically complete CAT(1) space.
\elem

\proof
Since $X$ contains a 2-flat, the essential Tits boundary is non-empty. 
Let $\xi$ and $\hat\xi$ be points in $\etits X$. Then there is a geodesic $\al$ of length at most $\pi$ between them in $\tits X$.
We will show that this geodesic is contained in $\etits X$. Since $\etits X$ is geodesically complete, we may assume that $\al$ has length $\pi$
to begin with. Let $\eta$ be a point on $\al$ and let $\eps>0$ be such that $\eta\notin N_\eps(\{\xi,\hat\xi\})$.
Now we extend $\al$ beyond $\xi$ to a point $\xi'$ by a geodesic $\al'$ of length $\eps$ in $\etits X$.
Let $\beta$ be a geodesic in $\tits X$ from $\xi'$ to $\hat \xi$. 
Then $\xi\notin\beta$ and  since $\tits X$ is CAT(1), we must have 
$\beta\cap(\al\setminus B_\eps(\hat\xi))=\emptyset$.
The union $\al\cup\al'\cup\beta$ contains a simple closed local geodesic $c$.
By definition, $c$ lies in $\etits X$, and by construction, $c$ contains $\al\setminus B_\eps(\hat\xi)$ and therefore $\eta$. 
Thus, all of $\al$ is contained in $\etits X$ as required. 

Now, the convexity implies that $\etits X$ is geodesic. Since $\etits X$ cannot contain a closed local geodesic of length less then $2\pi$
it is a CAT(1) space of dimension 1. By Lemma~\ref{lem_gc}, it is geodesically complete.
\qed

\section{Closing Lemma}\label{sec_cl}

\subsection{Asymptotically minimal actions}

\bdfn
We call an action $\Ga\acts X$ on a CAT(0) space  {\em asymptotically minimal}, if the induced action $\Ga\acts\geo X$ is {\em minimal}, i.e. 
$\geo X$ does not contain
closed $\Ga$-invariant proper subsets. 
\edfn

The following is our source for periodic rank 1 geodesics.

\bprop[{\cite[Proposition~1.10]{BaBu_periodic}}]\label{prop_BB}
Suppose $\Ga$ acts geometrically and asymptotically minimal on a CAT(0) space $X$. If the diameter of the Tits boundary $\tits X$ is strictly larger than
 $\pi$, the $X$ contains a $\Ga$-periodic geodesic of rank one.
\eprop

In order to deal with actions which are not asymptotically minimal, we will rely on work of Russell Ricks \cite{Ricks_1D}.

Recall the following construction due to Guralnik-Swenson \cite{GS_trans}.
Let $\Ga$ be a discrete group acting on
a compact Hausdorff space $Z$. Denote by $\beta\Ga$ the Stone–Čech compactification of $\Ga$. By compactness, for each $z\in Z$,
we can extend 
the orbit map $\rho_z:z\mapsto\ga z$ to a map $\beta\rho_z$ and, for $\om\in\beta\Ga$ we define 
\[T^\om:Z\to Z,\ z\mapsto \beta\rho_z(\om).\]
For fixed $z\in Z$, the
map $\beta G\to Z$ which maps $\om$ to $T^\om(z)$ is  continuous.
The family $\{T^\om\}_{\om\in\beta\Ga}$
of operators is closed under composition. 

%The inverse map g 󳨃→ g
%−1 on G extends to
%a continuous involution S : βG → βG; however, T
%ωT
%Sω usually only equals the identity for ω ∈ G.

Now let $\Ga\acts X$ be a geometric action on a locally compact CAT(0) space. 
Since $\bar X= X\cup\geo X$ is a compact Hausdorff space, the above construction applies.
Guralnik and Swenson observed that every operator $T^\om$ is 1-Lipschitz on $\tits X$ by the semi-continuity of the Tits metric.

Following \cite{Ricks_1D}, if $\si\subset\tits X$ is a round sphere and $\om\in\beta \Ga$ is such that  $T^\om:\tits X\to T^\om\si$
is a retraction which restricts to an isometry $\si\to T^\om \si$, we say $\om$ folds $\tits X$ (and $\si$) onto $T^\om\si$.
We call $T^\om\si$ a {\em folded round sphere}.
For every  top-dimensional round sphere $\si\subset\tits X$ there exists $\om\in\beta\Ga$ and a folded sphere $\si'$ such that $T^\om$
restricts to an isometry $\si\to\si'$ \cite[Corollary~3.2]{Ricks_1D}.

Let $X$ be a locally compact CAT(0) space with a geometric group action $\Ga\acts X$ and without $3$-flats.
Further, let $M\subset X$ be a minimal closed $\Ga$-invariant set. If $M$ intersects a round circle $\si\subset\tits X$ in an infinite set,
then the action $\Ga\acts\geo X$ is minimal, $M=\geo X$ \cite[Lemma~4.1]{Ricks_1D}.

\blem[{\cite[Lemma~4.8, Theorem~5.2]{Ricks_1D}}]\label{lem_ricks}
Let $\Ga\acts X$ be a geometric action on a locally compact CAT(0) space without $3$-flats.
Suppose that $\si\subset\tits X$ is a folded round circle.
Let $M\subset\geo X$ be a minimal closed $\Ga$-invariant proper subset which minimizes $l=\# (M\cap\si)$
among all such sets. If $l\geq 3$ is an odd number, then the following holds.
\begin{enumerate}
	\item $\etits X\subset\tits X$ is a closed convex subset.
	\item $\etits X$ is a union of round circles.
	\item Every round circle in $\etits X$ contains precisely $2l$ branch points,  evenly spaced around the circle;
	these branch points coincide with the points of $M\cup\Ant(M)$ on the circle.
	\item $\tits X$ is the union of $\etits X$ with a (possibly empty) collection of trees.
	Each such tree $T$ intersects $\etits X$ at a single point $m=m(T)\in M$ and $T\subset B_{\frac{\pi}{l}}(m)$.	
\end{enumerate}
If $l<3$ or $l$ even, then  $\tits X$ is a spherical building.
Furthermore, if $\diam(\tits X)=\pi$, then $\tits X$ is also a spherical building.
\elem

\subsection{Buildings at infinity}

It is well-known that the Tits boundary of a higher rank symmetric space or higher dimensional Euclidean building is a spherical building 
\cite[Proposition~4.2.1]{KleinerLeeb}.
Vice versa, Bernhard Leeb proved the following striking result.

\bthm[{\cite[Main Theorem]{Leeb}}]\label{thm_leeb}
Let $X$ be a locally compact, geodesically complete CAT(0) space. If $\tits X$ is a connected thick irreducible spherical building,
then $X$ is either a symmetric space or a Euclidean building.
\ethm

In the presence of a geometric group action one can relax the assumptions slightly, as we are now going to show.

\bprop\label{prop_buil}
Let $X$ be a locally compact CAT(0) space with a geometric group action $\Ga\acts X$. Further, let
$B\subset \tits X$ be a spherical building with $\dim(B)=\dim(\tits X)\geq 1$. 
If $\tits X\subset N_{\frac{\pi}{2}}(B)$, then $\tits X=B$.
\eprop
\medskip

%Before we turn to the proof, recall the following obvious fact.
%
%
%\blem\label{lem_boundclosed}
%Let $X$ be a locally compact CAT(0) space.
%Let $C\subset X$ be a closed convex subset. 
%Then $\tits C$ is closed in $\geo X$. 
%\elem
%
%\proof
%Immediate since $X$ is locally compact.
%\qed
%\medskip

The proof of Proposition~\ref{prop_buil} relies on the following observation.

\blem\label{lem_closed}
Let $X$ be a locally compact CAT(0) space and let
$B\subset \tits X$ be a spherical building with $\dim(B)=\dim(\tits X)\geq 1$ and $\tits X\subset N_{\frac{\pi}{2}}(B)$. 
Then $B$ is closed in $\geo X$.
\elem

\proof
Let $(x_k)$ be a sequence in $B$ which converges to a point $\xi$ in $\geo X$.
Since $B$ is top-dimensional, $\xi$ has an antipode $\hat\xi$ in $B$.
We extend geodesics $\hat\xi x_k$ inside $B$ to antipodes $\xi_k$ of $\hat\xi$.
Choose geodesics $c_k\subset B$ from $\hat\xi$ to $\xi_k$ which all have the same starting direction and therefore agree up to a certain point $\eta\in B$.
By compactness of $\geo X$ and semi-continuity of the Tits metric, we can pass to a subsequence such that $c_k$ converges to a geodesic $c$ from $\hat\xi$
to $\xi$. Note that $c$ contains the segment $\hat\xi\eta$. 
Now extend the segment $\eta\hat\xi$ up to a point $\hat\zeta\in B$ such that 
$\hat\xi$ and $\hat\eta$ lie in the same chamber $\hat\Delta\subset B$.
Denote by $\zeta_k$ the antipode of $\hat\zeta$ on $c_k$. Then $\xi_k$ and $\zeta_k$ lie in a chamber $\Delta_k\subset B$. 
Let $\al_k\subset B$ be an apartment which contains the chambers $\hat\Delta$ and $\Delta_k$.
Then there is a round circle $\si_k\subset \al_k$
which contains the four points  $\hat\zeta,\hat\xi,\zeta_k,\xi_k$.
Moreover, since $|\hat\xi,\zeta_k|<\pi$, we also have $\eta\in\si_k$.
Let $f_k:S^1\to\tits X$ be an isometric embedding with image $\si_k$.  Using compactness of $\geo X$ and semi-continuity of the Tits metric again, 
we obtain a 1-Lipschitz limit map $f_\infty:S^1\to\tits X$. Since $\hat\xi\hat\zeta\subset\si_k$, $\hat\xi \eta\subset\si_k$  and $\xi_k\to\xi$, we see that the image $\si_\infty$ of $f_\infty$ contains the geodesic $c$ and the segment $\hat\xi\hat\zeta$. It follows that $\si_\infty$
is a round circle. Since $\tits X\subset N_{\frac{\pi}{2}}(B)$, $\si_\infty$ has to lie in $B$. Thus $\xi\in B$ as required.
\qed

\begin{rem}
Note that the strict $\frac{\pi}{2}$-density above is necessary.
\end{rem}

\proof[Proof of Proposition~\ref{prop_buil}]
Every round sphere $\si\subset\tits X$ intersects $B$ in a convex set $C$.
By assumption, the set $C$ cannot be contained in a hemisphere of $\si$.
Thus $\si$ lies entirely in $B$. If $\si$ is top-dimensional, then, by \cite[Corollary 2.5]{BMS_affine}, the orbit $\Ga\si$
is dense in $\geo X$ \footnote{\cite[Corollary 2.5]{BMS_affine} only states that the family of top-dimensional round spheres
is dense in $\geo X$ but their argument proves the stronger form.}. Since $B$ is closed by Lemma~\ref{lem_closed}, the conclusion follows.
\qed

\subsection{Diameter bound}

By Proposition~\ref{prop_BB}, in order to prove the Closing Lemma it is enough to verify the following.

\bthm\label{thm_diambound}
Let $X$ be a locally compact CAT(0) space without 3-flats. Suppose that
$\Ga\acts X$ is a geometric action which is not asymptotically minimal.
Then the Tits boundary $\tits X$ is a spherical join or a spherical building. 
In particular, the diameter of $\tits X$ is equal to $\pi$.
\ethm

For the rest of this section, the assumptions of Theorem~\ref{thm_diambound}
are in place. Let $\si\subset \tits X$ be a folded round circle.
We choose  a minimal closed $\Ga$-invariant proper subset $M\subset\geo X$  which minimizes $l=\# (M\cap\si)$
among all such sets. Recall that the number $l$ has to be finite, otherwise $\Ga$ acts minimally \cite[Lemma~4.1]{Ricks_1D}.
Moreover, we may assume $l$ is odd and $\geq 3$ \cite[Theorem~D]{Ricks_1D}.
Recall that the {\em essential Tits boundary} $\etits X$ of $X$ is the subset of $\tits X$
given by the union of all simple closed local geodesics.

\blem\label{lem_esstits}
The essential Tits boundary $\etits X$ is an irreducible spherical building.
\elem

\proof
We show that the diameter $D$ of $\etits X$ is equal to $\pi$.
By Lemma~\ref{lem_ricks}, $\etits X$ has the structure of a simplicial complex with edge length $\frac{\pi}{l}$ and diameter at most $\frac{l+1}{l}\pi$.
Moreover, every edge has precisely one boundary vertex in $M$
and one boundary vertex in $\Ant(M)$.  Thus, if $D$  is strictly larger than $\pi$,
then $M\cap\Ant(M)\neq\emptyset$. But then $M\subset\tits X$ is a proper closed subset with $\Ant(M)\subset M$ \cite[Lemma~3.17]{Ricks_1D}.
In turn, $\etits X$ would have to be a spherical join or spherical building \cite{Ly_rigidity}. Contradiction.
Thus $\etits X$ has diameter $\pi$. The claim follows from \cite[Theorem~6.1]{CL_metric}.
\qed

\proof[Proof of Theorem~\ref{thm_diambound}]
By Lemma~\ref{lem_esstits}, $B:=\etits X$ is an irreducible spherical building. 
By Lemma~\ref{lem_ricks}, $B$ is $\frac{\pi}{3}$-dense in $\tits X$.
Thus, $B\subset\geo X$ is closed by Lemma~\ref{lem_closed}. Therefore
we conclude $\tits X=B$ from Proposition~\ref{prop_buil}.  
\qed

\section{Diameter Rigidity}\label{sec_rank_2}

\subsection{Projecting ideal points onto parallel sets}

\blem\label{lem_minimum}
Let $X$ be a locally compact CAT(0) space with $1$-dimensional Tits boundary. 
Let $\xi^-$ and $\xi^+$
be a pair of antipodes in $\tits X$ and let $\eta$ be a pole of $\{\xi^-,\xi^+\}$.
Further, let $P\subset X$ be a closed convex subset with $\eta\notin\tits P$ and $\{\xi^-,\xi^+\}\subset\tits P$.
If $b_\eta$ is a horofunction associated to $\eta$, then $b_\eta$ attains a minimum on $P$. Moreover, the minimum set 
contains a complete geodesic $c$ with $\geo c=\{\xi^+,\xi^-\}$.  
\elem

\proof
Let us first show that the parallel set $P(\xi^+,\xi^-)$ is non-empty. Denote by $\tau^+\subset \tits X$
the round hemisphere with $\D\tau^+=\{\xi^+,\xi^-\}$ and center $\eta$.
If there is no round hemisphere $\tilde\tau^+$ in $\tits P$ with $\D\tilde\tau^+=\{\xi^+,\xi^-\}$, then by \cite[Proposition~2.1]{Leeb},
there exists a complete geodesic $c$ in $P$ with $\geo c=\{\xi^+,\xi^-\}$. On the other hand, if there does exist such a 
$\tilde\tau^+$, then, since  $\eta\notin\tits P$, the union $\si:=\tilde\tau^+\cup\tau^+$ forms a round 1-sphere in $\tits X$.
Since $\dim(\tits X)=1$, \cite[Proposition~2.1]{Leeb} ensures the existence of a 2-flat $F$ with $\tits F=\si$ and therefore $F\subset P(\xi^+,\xi^-)$.
Now it follows from \cite[Sublemma~2.3]{Leeb}, if $HB_\eta$ is a horoball based at $\eta$ such that $HB_\eta\cap P$ is non-empty, then 
$HB_\eta\cap P$ contains a complete geodesic $c$  with $\geo c=\{\xi^+,\xi^-\}$.
Since $b_\eta$ is constant on such geodesics, in order to find the required minimum, it is enough to restrict $b_\eta$
to $P\cap CS(\xi^+,\xi^-)$. Note that $|\eta',\eta|=\pi$ for every point $\eta'\in\tits CS(\xi^+,\xi^-)\setminus\{\eta\}$. 
Hence $\lim\limits_{x\to\eta'}b_\eta(x)=+\infty$ and 
$b_\eta$ does attain a minimum at a point $p$ in $P$. 
By assumption, we have
\[\{\xi^-,\xi^+\}\subset\tits(HB(\eta,p)\cap P).\]
Thus, by \cite[Sublemma~2.3]{Leeb}, there is a  complete geodesic $c$ in $HB(\eta,p)\cap P$ with $\geo c=\{\xi^+,\xi^-\}$.
Since $b_\eta$ is bounded above on $HB(\eta,p)$ by its value at $p$  the claim follows.   
\qed

\blem\label{lem_proj}
Let $X$ be a locally compact CAT(0) space with $1$-dimensional Tits boundary. Let $\si\subset\tits X$ be a round 1-sphere.
Let $\tau^+\subset\tits X$ be a round hemisphere with $\tau^+\cap\si=\D\tau^+$.
Then there exists a 2-flat $F\subset X$ with $\tits F=\si$ and a flat half-plane $H$,  orthogonal to $F$ and with
$\tits H=\tau^+$. Moreover, for any geodesic $c\subset F$ which is not parallel to $\D H$, holds $\angle(H,P(c))\geq\frac{\pi}{2}$.
\elem

\proof
Since $\dim(\tits X)=1$, the parallel set $P(\si)$  is a non-empty closed convex subset of $X$
which splits isometrically as $P(\si)\cong\R^2\times C$ where $C$ is a compact CAT(0) space.
Let $\eta$ be the center of $\tau^+$ and let $b_\eta$ be an associated horofunction.
By Lemma~\ref{lem_minimum}, 
$b_\eta$ attains a minimum on $P(\si)$ and the minimum set $Z\subset P(\si)$ contains a complete geodesic $l$ with $\geo l=\D\tau^+$.
Let $H$ be the flat half-plane with $\tits H=\tau^+$ and $\D H=l$.
Because $l\subset P(\si)$, there exists a 2-flat $F\subset X$ with $\tits F=\si$ which contains $l$.
Since $l$ lies in the minimum set $Z$, the half-plane $H$ has to be orthogonal to $F$.

Let $c\subset F$ be a geodesic not parallel to $l$. Choose  a point $x\in l$ and a geodesic ray $\rho$ starting in $x$ and asymptotic to $\eta$.
In particular, $\rho\subset H$.
 
We first claim that $\rho(t)\notin P(c)$ for $t>0$. Indeed, if $\rho(t)\in P(c)$, then, since $\rho(t)\in P(l)$, there exists lines $c_t\| c$ and $l_t\| l$ through
$\rho(t)$. These span a flat $F_t\subset P(\si)$. Hence $t=0$ because $\rho(0)$ minimizes $b_\eta$.

Now suppose $\angle(\rho,P(c))<\frac{\pi}{2}$. Then there exists a geodesic $c_1$ in $P(c)$ starting in $x$ and realizing the angle, 
$\angle(\rho,c_1)=\angle(\rho,P(c))$. In particular, $b_\eta(c_1(t))<b_\eta(x)$. 
Note that since $c$ is not parallel to $l$, the point $\eta$ is not in $\tits P(c)$. Moreover, we have $\D\tau^+\subset\tits P(c)\cap\bar B_{\frac{\pi}{2}}(\eta)$.
Hence, by Lemma~\ref{lem_minimum}, $b_\eta$ attains a minimum on $P(c)$, and the minimum set $Z_c$ contains a complete geodesic 
$\tilde l$ with $\D \tilde l=\{\xi, \hat\xi\}$.
But then $\tilde l$ lies in $P(c)$ and therefore, arguing as above, there is a 2-flat $\tilde F\subset P(\si)$ with $\tilde l\subset \tilde F$.
This contradicts the fact that $Z$ minimizes $b_\eta$ on $P(\si)$. Therefore $\angle(H,P(c))=\angle(\rho,P(c))\geq\frac{\pi}{2}$.
\qed

\subsection{Parallel sets do not accumulate in the Tits metric}\label{subsec_no_acc}

\blem\label{lem_no_acc}
Let $X$ be a locally compact CAT(0) space with $1$-dimensional Tits boundary. 
Let $g$ be  an axis of an axial isometry $\ga$. 
Suppose that there exists $a>0$ and a sequence of 2-flats $(F_k)$ in $P(g)$
with $g\subset N_a(F_k)$ for every $k\in\N$. Further, suppose that there is a sequence $(H_k)$
of flat half-planes with boundaries $h_k:=\D H_k\subset F_k$ and $\geo H_k\cap\geo F_k=\geo h_k$.
If $\geo h_k\to\geo g$ with respect to the Tits metric, then $\geo h_k=\geo g$
for almost all $k$.
\elem

\proof
Suppose for contradiction that the convergence $\geo h_k\to\geo g$ is non-trivial.
By Lemma~\ref{lem_proj}, there exist 2-flats $\tilde F_k\subset P(g)$ parallel to $F_k$; 
and flat half-planes $\tilde H_k$ orthogonal to $\tilde F_k$ and with $\geo \tilde H_k=\geo H_k$; and such that 
 $\tilde H_k$ is orthogonal to $P(g)$. Set $\tilde h_k:=\D \tilde H_k$. Since $X$ has rank 2 and a cocompact isometry group, there exists $D>0$
such that $|F_k,\tilde F_k|_H\leq D$ for all $k\in\N$. In particular, 
there exist complete geodesics $\tilde g_k\subset  \tilde F_k$ which are parallel to $g$ and
satisfy $|g, \tilde g_k|\leq a+D$ for all $k\in\N$. Denote by $\tilde x_k\in \tilde F_k$ the intersection point of $\tilde g_k$ and $ \tilde h_k$.
Then there are powers of $\ga$, denoted by $\ga_k$, such that the points $\ga_k(\tilde x_k)$ lie in a fixed compact set.
Since $\ga$ preserves $P(g)$, $\ga_k \tilde H_k$ is still orthogonal to $P(g)$. Since $X$ is locally compact, after choosing a subsequence, 
the flat half-planes $\ga_k \tilde H_k$ converge to a flat half-plane $\tilde H_\infty$. 
Note that by assumption $\ga_k (\geo \tilde h_k)\to\geo g$ with respect to the
Tits metric.
Hence $\tilde H_\infty\subset P(g)$. On the other hand, by upper semi-continuity of angles, $\tilde H_\infty$ has to be orthogonal to $P(g)$. Contradiction.
\qed

\subsection{Regular points at infinity}

In this section we will prove that if the Tits boundary of a cocompact CAT(0) space has dimension $1$ and diameter $\pi$, then it
contains regular points.
To achieve this, we will use the following technical result in an essential way.

\begin{namedlemma}[Half-Plane]\label{lem_key_tech}
Let $X$ be a locally compact CAT(0) space with $1$-dimensional Tits boundary.  Let $F\subset X$ be a 2-flat with $\geo F=\si$.
Further, let $g$ be an axis of an axial isometry $\ga$ and assume $\geo g\subset\si$.
Let $\xi^\pm\in\si$ be antipodes disjoint from $\geo g$.
Suppose that there is a sequence of local geodesics $\al^+_k\subset \tits X$ such that 
\begin{itemize}
	\item $\D\al^+_k=\{\xi^-_k,\xi^+_k\}$ and $\xi^\pm_k\to\xi^\pm$ with respect to the Tits metric;
	\item the length of $\al_k^+$  converges to $\pi$ as $k\to\infty$.
\end{itemize}
Then there exists a 2-flat $\tilde F\subset P(g)$, and a flat half-plane $\tilde H$ with boundary $\tilde h:=\D \tilde H\subset\tilde F$ such that the following properties hold.
\begin{enumerate}
	\item $g\subset N_a(\tilde F)$ if $g\subset N_a(F)$ for $a>0$;
	\item $|\geo \tilde h,\geo g|=|\{\xi^-,\xi^+\},\geo g|$;
	\item $\geo\tilde H\cap\geo\tilde F=\geo\tilde h$.
\end{enumerate}
\end{namedlemma}

The  proof of the \hyperref[lem_key_tech]{Half-Plane Lemma} requires methods from geometric measure theory.
Since these techniques do not play a role in the rest of the paper,
we defer their discussion, as well as the proof of the \hyperref[lem_key_tech]{Half-Plane Lemma}, to Appendix~\ref{sec_app}. 
\medskip

The following is certainly well-known, we include it for completeness.

\blem\label{lem_per2fl}
Let $X$ be a locally compact CAT(0) space with a geometric group action $\Ga\acts X$. 
Let $F\subset X$ be a $\Ga$-periodic $k$-flat which bounds a flat $(k+1)$-half-space $H$.
Suppose that $G\cong\Z^k$ is a subgroup in $\Ga$ which acts geometrically on $F$.
If $G$ preserves $H$, then there exists a $\Ga$-periodic $(k+1)$-flat $\hat F\subset P(F)$.
\elem

\proof
Let $\rho\subset H$ be a geodesic ray orthogonal to $F$. Choose points $x_k\in\rho$ with $|x_k,F|\to\infty$.
We find a sequence $(\ga_k)\subset\Ga$ such that $\ga_k(x_k)$ stays in a fixed compact set $K$.
We can pass to a subsequence such that for all $k<l$ the element $\ga_l^{-1}\ga_k$
is axial (cf. proof of \cite[Theorem~11]{Sw_cut}). Now note that since $G$ preserves $H$, each of its elements $g$ has  constant   
displacement on $H$, $|x,g x|=a_g>0$ for every $x\in H$.
Thus $|\ga_k x_k,(\ga_k g\ga_k^{-1})\ga_k x_k|=|x_k, g x_k|=a_g$. Since $\Ga$ acts properly discontinuously and $\ga_k x_k\in K$, we may pass to a further subsequence
 such that $\ga_1 g\ga_1^{-1}=\ga_k g\ga_k^{-1}$ holds for all $k\in \N$.
Hence $\beta_k=\ga_k^{-1}\ga_1$ is an axial isometry which commutes with $g$. Letting $g$ vary through a basis of $G$ implies the claim.
% The claim follows from \cite[Flat Torus Theorem~7.1]{BH}.
\qed

\bprop\label{prop_regular_points}
Let $X$ be a locally compact CAT(0) space with $1$-dimensional Tits boundary and a geometric group action $\Ga\acts X$. 
Suppose that the diameter of $\tits X$ is equal to $\pi$.
Let $g$ be an axis of an axial isometry $\ga\in\Ga$. Then there exists a positive $\eps$ and a round 1-sphere $\si\subset\tits P(g)$ 
such that  $\dot B_\eps(g(+\infty))\cap\si$ 
does not contain a branch point. In particular, $\si\subset\tits X$ contains a non-empty open relatively compact subset.
\eprop

\proof
Since the diameter of $\tits X$ is $\pi$, the axis $g$ bounds a flat half-plane $H\subset X$.
If $\ga$ preserves $H$, then by Lemma~\ref{lem_per2fl}, the parallel set $P(g)$ contains a periodic 2-flat. 
In this case, the claim follows from  \cite[Corollary~5.7]{St_rrI}.
If $H$ is not preserved, then $\si=\tau^+\cup\ga\tau^+$ is a round sphere in $\tits X$ where $\tau^+:=\geo H$.
Let $F$ be a 2-flat with $\geo F=\si$. Let $a>0$ be such that $g\subset N_a(F)$.
Now suppose that there is a sequence $(\xi_k^+)$ of pairwise distinct branch points in $\si$ with $\xi_k^+\to g(+\infty)$.
Denote by $\xi_k^-$ the antipode of $\xi_k^+$ in $\si$.
Since $\tits X$ has diameter $\pi$, we find for each $k\in\N$ a sequence of local geodesics $\al_{kl}\subset \tits X$ such that 
\begin{itemize}
	\item $\D\al_{kl}=\{\xi^-_{kl},\xi^+_{kl}\}$ and $\xi^\pm_{kl}\to\xi_k^\pm$ with respect to the Tits metric;
	\item the lengths of $\al_{kl}$  converges to $\pi$ as $l\to\infty$.
\end{itemize}
Hence by the \hyperref[lem_key_tech]{Half-Plane Lemma}, there exists a sequence of 2-flats $(\tilde F_{k})$ in $P(g)$, 
and a sequence of flat half-planes $(\tilde H_{k})$ with boundaries $\tilde h_k:=\D\tilde H_k\subset\tilde F_k$ and such that
the following properties hold.
\begin{enumerate}
	\item $g\subset N_a(\tilde F_k)$ for every $k\in\N$;
	\item $|\geo \tilde h_k,\geo g|=|\{\xi_k^-,\xi_k^+\},\geo g|$;
	\item $\geo\tilde H_k\cap\geo\tilde F_k=\geo\tilde h_k$.
\end{enumerate}
In particular, $\geo \tilde h_k\to \geo g$ with respect to the Tits metric.
Hence Lemma~\ref{lem_no_acc} implies that $\tilde h_k$ is parallel to $g$ for almost all $k$.
This is a contradiction since the points $\xi_k^+$ are pairwise distinct.
\qed

\subsection{The Tits boundary is a building}

\blem\label{lem_reg_1D}
Let $Z$ be a 1-dimensional CAT(1) space of diameter $\pi$.
Then the subset $O\subset Z$ of regular points is closed under taking antipodes.
\elem

\proof
Let $\xi\in O$ be a regular point and let $\hat\xi\in Z$ be an antipode.
Then there exists $s>0$ such that $\bar B_s(\xi)$ is isometric to a closed interval.
Denote by $\xi^\pm$ the two endpoints of $\bar B_s(\xi)$.
Since the diameter of $Z$ is equal to $\pi$, we see that
$\xi$ is the unique antipode of $\hat\xi$ in $B_s(\xi)$. 
Thus, $|\xi,\xi^\pm|+|\xi^\pm,\hat\xi|=\pi$ and $\xi,\hat\xi$ lie in a round 1-sphere $\si$.
We claim that $B_s(\hat\xi)\subset\si$. Let $\hat\eta\in Z$ be  a point in  $B_s(\hat\xi)$.
Then $\hat\eta$ has an antipode $\eta\in B_s(\xi)$ and as before, it has to be unique and we find a round 1-sphere $\si'$
which contains the points $\eta$ and $\hat\eta$. By regularity, $\bar B_s(\xi)\subset\si\cap\si'$.
Since $|\hat\xi,\eta|<s$, we conclude $\si=\si'$ and therefore $\eta\in\si$ as required.
\qed

\blem\label{lem_etits_build}
Let $X$ be a locally compact CAT(0) space with $1$-dimensional Tits boundary and a geometric group action. 
Suppose that the diameter of $\tits X$ is equal to $\pi$. Then $\etits X$ is a spherical join or a spherical building.
\elem

\proof
We may assume that $\etits X$ is not a round sphere. By \cite[Theorem~C]{Kleiner}, $X$ contains a 2-flat.
Thus, by  Lemma~\ref{lem_diameter_etits}, the diameter of $\etits X$ is equal to $\pi$. By \cite[Theorem~11]{Sw_cut}, the group $\Ga$
contains an axial element.
Hence, the open subset $O\subset\etits X$ of regular points is non-empty by Proposition~\ref{prop_regular_points}.
We infer from Lemma~\ref{lem_reg_1D} that $\etits X\setminus O$ is a proper closed subset of $\etits X$ which contains with every point all of its antipodes.
The claim follows from \cite[Main Theorem]{Ly_rigidity}.
\qed

\bcor\label{cor_rank2_build}
Let $X$ be a locally compact CAT(0) space with $1$-dimensional Tits boundary and a geometric group action. 
Suppose that the diameter of $\tits X$ is equal to $\pi$. Then $\tits X$ is a spherical join or a spherical building.
\ecor

\proof
By Proposition~\ref{prop_buil} and Lemma~\ref{lem_etits_build}, it is enough to show that every point in $\tits X$ has distance less than $\frac{\pi}{2}$
from $\etits X$. Suppose for contradiction that there exists a point  $\xi\in\tits X\setminus\etits X$ at 
distance $\frac{\pi}{2}$ from $\etits X$. Choose $\eta\in\etits X$  such that 
$|\xi,\eta|<\frac{\pi}{2}+\delta$ for a small positive $\delta$.
Extend the geodesic $\xi\eta$ inside $\etits X$ up to a local geodesic $\al$ of length $\pi+\delta$ and let $\eta'$ denote its endpoint.
Now choose a geodesic $\beta$ from $\eta'$ to $\xi$. Since $\tits X$ is CAT(1), $\beta$ avoids the interior of  a subsegment $\al_0\subset\al$
of length $\pi$. In particular, $\al\cup\beta$ contains a simple closed local geodesic $c$ of length at least  $2\pi$.
Thus, the geodesic $\xi\eta$ intersects $\etits X$ in a segment of length at least $\frac{\pi}{2}-\delta$ and distance between $\xi$
and $\etits X$ is less then $\frac{\pi}{2}$. Contradiction.
\qed

\proof[Proof of Theorem~\ref{thm_mainB}]
By Corollary~\ref{cor_rank2_build}, the Tits boundary $\tits X$ is a spherical join or a building.
Since $X$ is geodesically complete and locally compact, the claim follows from Theorem~\ref{thm_leeb} and \cite[Proposition~2.3.7]{KleinerLeeb}.
\qed

\section{Applications}

In this section we provide the proofs of Corollaries~\ref{cor_tits} and~\ref{cor_main}.

\proof[Proof of Corollary~\ref{cor_tits}]
If the Tits diameter is strictly larger than $\pi$, then by Theorem~\ref{thm_mainA}, $X$ contains a $\Ga$-periodic geodesic of rank one.
Thus $\Ga$ contains a non-abelian free subgroup \cite[Theorem~3.5]{ballmannbook}. On the other hand, if the Tits diameter is equal to $\pi$,
then $X$ is either an irreducible symmetric space, an irreducible Euclidean building or splits metrically as a product.
The first case follows from Tits original theorem \cite{Ti_free} and the second case is covered by \cite[Theorem~F]{BaBr_orbi}.
We are left with the product case, $X\cong X_1\times X_2$. By assumption, the factors $X_i$ do not contain $2$-flats.
After passing to an index two subgroup, we may assume that $\Ga$ preserves the factors.
We will argue that there exist two $\Ga$-axes such that their projections to one of the factors
are not asymptotic to one another and therefore lead to a non-abelian free subgroup \cite[Theorem~3.5]{ballmannbook}.
Suppose for contradiction that any pair of $\Ga$-axes project to (one-sided) asymptotic geodesics in both factors.
Then every axis has to be regular, i.e. not parallel to a factor.
Since $\Ga$ acts geometrically, this means that every axial element lies in a subgroup $\Z^2<\Ga$.
Every such subgroup preserves a $2$-flat. Two different $\Ga$-periodic $2$-flats cannot be asymptotic to the same Weyl chamber because
$\Ga$ acts properly discontinuous. Thus, we either find the required $\Ga$-axes and hence the non-abelian subgroup, or $\Ga$ is
a Bieberbach group.
\qed

\proof[Proof of Corollary~\ref{cor_main}]
Under the assumptions, we find a $\Ga$-equivariant quasi-isometry $\Phi:X\to X'$.
Hence for every axial element $\ga\in\Ga$ a pair of $\ga$-axes $c\subset X$ and $c'\subset X'$
has to have the same rank. Indeed, $\Phi(c)$ lies at finite Hausdorff distance from $c'$
and therefore $c'$ has quadratic growth if $c$ does. Thus the claim follows from Theorems~\ref{thm_mainA} and~\ref{thm_mainB}
together with \cite[Theorem~1.1.3]{KleinerLeeb}.  
\qed

\appendix
\section{The relative asymptotic Plateau problem after Kleiner--Lang}\label{sec_app}

The primary goal of this section is to provide a proof of the \hyperref[lem_key_tech]{Half-Plane Lemma}.
Our proof requires a solution of an asymptotic Plateau problem relative to a flat.
Since this has not been done in the literature, we need some preparation.
The (non-relative) asymptotic Plateau problem has recently been solved  by Kleiner--Lang~\cite{KL_higher}
and we use their work extensively throughout this section. The idea to use solutions to a relative Plateau problem originated from \cite{HKS_I}.

In Section~\ref{subsec_notation}, we agree on notation and proof  monotonicity and volume rigidity of minimizing currents in CAT(0) spaces,
a result which is folklore.
In Section~\ref{subsec_relative}, we solve the asymptotic Plateau problem relative to a 2-flat.
We provide basic properties of relative minimizers which are required for our proof of the \hyperref[lem_key_tech]{Half-Plane Lemma}.
At last, in Section~\ref{subsec_rel_asym}, we present the proof of the \hyperref[lem_key_tech]{Half-Plane Lemma}.
We have refrained from any kind of general treatment, in particular we only discuss spaces of rank 2. 
This is all we need here and a broader discussion would get out of hand. 

The general reference for this section is \cite{KL_higher}.
See also \cite{AK, La_currents} for a thorough background on geometric measure theory in metric spaces.

\subsection{Currents in CAT(0) spaces}\label{subsec_notation}

Let $X$ be a locally compact metric space.
For every integer $n \ge 0$, let $\cD^n(X)$ denote the set of all 
$(n+1)$-tuples $(\pi_0,\ldots,\pi_n)$ of real valued functions on $X$ such 
that $\pi_0$ is Lipschitz with compact support $\spt(\pi_0)$ and 
$\pi_1,\dots,\pi_n$ are locally Lipschitz.
A {\em $n$-current} $S$ in $X$ is a function 
$S:\cD^n(X) \to \R$ satisfying the following three conditions:
\begin{enumerate}
	\item(multilinearity)
$S$ is $(n+1)$-linear;
\item(continuity) 
$S(\pi_{0,k},\ldots,\pi_{n,k}) \to S(\pi_0,\ldots,\pi_n)$
whenever $\pi_{i,k} \to \pi_i$ pointwise on $X$, $\bigcup_k\spt(\pi_{0,k})$ is bounded 
and the $\pi_{i,k}$ are uniformly locally Lipschitz continuous; 
\item(locality)
$S(\pi_0,\ldots,\pi_n) = 0$ whenever one of the functions
$\pi_1,\ldots,\pi_n$ is constant on a neighborhood of $\spt(\pi_0)$.
\end{enumerate}

We write $\cD_n(X)$ for the vector space of $n$-currents 
in $X$. For every $S \in \cD_n(X)$ we denote by  $\spt(S)$ the {\em support} of $S$ and by  $\D S \in \cD_{n-1}(X)$ its {\em boundary}.
As usual, $\|S\|$ will denote the associated regular Borel measure and $\M(S)=\|S\|(X)$ its {\em mass}. 
The {\em flat norm} of $S$ is denoted $\cF(S)$.
Moreover, for a Borel set $A\subset X$ we denote by $S\on A \in \cD_n(X)$ the {\em restriction} of $S$ to $A$.
For a proper Lipschitz map $f: X \to Y$ into another locally compact
metric space $Y$, we have the {\em push-forward} $f_\#S \in \cD_n(Y)$.
We denote by $\bI_{n,loc}(X)$ the abelian group of {\em locally integral $n$-currents},
and write $\bI_{n,c}(X)$ for the subgroup of {\em integral $n$-currents with compact support}.
Further, $\bZ_{n,loc}(X)\subset\bI_{n,loc}(X)$ and $\bZ_{n,c}(X)\subset\bI_{n,c}(X)$
will denote the subgroup of {\em cycles}, i.e. integral currents with boundary  zero.

Recall that an integral current $S\in \bI_{n,c}(X)$ is  concentrated on a countable $\Ha^n$-rectifiable set $\chi_S$, its {\em characteristic set} 
\cite[Theorem~4.6]{AK}. The characteristic set $\chi_S$ is unique in the sense that any Borel set $E$ with $\|S\|=\|S\|\on E$
contains $\chi_S$ up to $\Ha^n$-negligible sets. Moreover, there exists a $\Ha^n$-integrable {\em multiplicity function} $\theta:\chi_S\to\N$
such that $\|S\|=\theta\cdot\Ha^n\on\chi_S$ \cite[Theorem~9.5]{AK}. Note that by \cite[Lemma~9.2]{AK} the area factor $\la$ appearing in the general formula in 
\cite[Theorem~9.5]{AK} is equal to $1$ in our case, since CAT(0) spaces have Euclidean tangent cones, cf.~\cite[Section~11]{LW_plateau}.

On  Euclidean space $\R^m$ we denote the {\em flat $n$-chains with compact support} 
by $\cF_{n,c}(\R^m)$, see \cite{Fl_flat, W_defo, W_rec}.
Recall that for a compact subset $K\subset\R^m$ the space of flat chains with compact support in $K$ is the flat-closure
of $\cP_{n}(K)$, {\em the polyhedral $n$-chains with support in $K$}.
We will need a localized version. Define
 a {\em local flat $n$-chain} as an $n$-current $S$
such that for every $x\in\R^m$ there exists $S'\in\cF_{n,c}(\R^m)$ with $x\notin\spt(S-S')$.
The space of  {\em locally flat $n$-chains} is denoted by $\cF_{n,loc}(\R^m)$.

Let  $X$ be a CAT(0) space and $S\in \bZ_{n,c}(X)$. Then for every point $p\in X$ we have the {\em cone}
from $p$ over $S$, denoted by $C_p S\in\bZ_{n+1,c}(X)$, cf. \cite[Section~2.7]{KL_higher}. If $\spt(S)\subset \bar B_R(p)$ the following
{\em Euclidean cone inequality} holds.
\[\M(C_p S)\leq \frac{R}{n+1}\cdot\M(S).\]  

The fundamental class of $\R^n$ is denoted by $\bb{\R^n}\in\bZ_{n,loc}(\R^n)$. Every proper Lipschitz map $\varphi:\R^n\to X$
induces a natural local cycle $\varphi_\#(\bb{\R^n})\in\bZ_{n,loc}(X)$. In this way, a flat $F$ in $X$ becomes a local cycle which we will still denote by $F$
and call a {\em multiplicity 1 flat}.

Following \cite{KL_higher}, for a current
$S \in \bI_{n,\loc}(X)$, a point $p \in X$ and $r > 0$, we define the {\em ($r$-)density at $p$} by 
\[
\G_{p,r}(S) := \frac{1}{r^n} \|S\|(B_r(p))\quad \G_{p}(S) :=\liminf\limits_{r\to 0} \frac{1}{r^n} \|S\|(B_r(p)).
\]  
 Furthermore, for any $p \in X$, we define the  {\em density at infinity} by  
\[
\Gi(S) := \limsup\limits_{r\to\infty} \G_{p,r}(S).
\]
Similarly, we define the {\em filling density at $p$} by 
\[
\F_{p,r}(S) := \frac{1}{r^{n+1}} 
\inf\{\M(V) : V \in \bI_{n+1,c}(X),\,\spt(S-\D V) \cap B_r(p) =\emptyset\}.
\]  
Furthermore, for any $p \in X$, we define the  {\em filling density at infinity} by 
\[
\Fi(S) := \limsup\limits_{r\to\infty} \F_{p,r}(S).
\]

The following monotonicity property of minimizers is certainly well-known. We include
a version for CAT(0) spaces, since we were unable to find a reference. 

\bprop[monotonicity]\label{prop_mon}
Let $X$ be a locally compact CAT(0) space and $S \in \bI_{n,c}(X)$ an area minimizer.
Then for every point $p\in\spt(S)\setminus\spt(\D S)$  the $r$-density $\G_{p,r}(S)$
is a non-decreasing function of $r$, as long as $r\leq|p,\spt(\D S)|$, and bounded below by $\om_n$,
\[\om_n\leq \G_{p,r}(S)\nearrow.\]
Moreover, if $\G_{p,r_0}(S)=\om_n$ holds for some $r_0\leq|p,\spt(\D S)|$, then $\spt(S)\cap \bar B_{r_0}(p)$ is 
isometric to a Euclidean ball of radius $r_0$ and $S$ has constant multiplicity 1.
\eprop

\proof
Monotonicity follows from the standard cone comparison since CAT(0) spaces satisfy a Euclidean coning inequality.
More precisely, for almost all $r\leq |p,\spt(\D S)|$ holds 
\[\M(S\on B_r(p))\leq \frac{r}{n}\cdot\M(\D(S\on B_r(p))).\]
By the coarea inequality we have 
\[\frac{d}{dr}\M(S\on B_r(p))\geq \M(\D(S\on B_r(p)))\] 
for almost all $r$, and monotonicity follows by integration.
Recall that $S$ is concentrated on the countably $\Ha^n$-rectifiable characteristic set $\chi_S$ which is given by $\chi_S:=\{x\in X|\ \G_x(S)>0\}$.
The isoperimetric inequality yields some positive lower density bound at all points $p\in\spt(S)\setminus\spt(\D S)$ \cite[Lemma~3.3]{KL_higher}.
This implies $\chi_S\cap\bar B_r(p)=\spt(S)\cap\bar B_r(p)$ for every $r<|p,\spt(\D S)|$.
It follows from \cite[Theorem~5.4]{AK_rec} that $\G_x(S)=\theta\cdot\om_n$ for almost all points $x\in\chi_S$ where $\theta:\chi_S\to \N$ denotes the multiplicity function of $S$.
The bound  $\G_x(S)\geq\om_n$ extends to the set $\spt(S)\setminus\spt(\D S)$ by upper semi-continuity of $\G_x(S)$.

Now suppose $\G_{p,r_0}(S)=\om_n$ holds for some $r_0<|p,\spt(\D S)|$.
Then the above inequalities become equalities and we have
\[\M(S\on B_r(p))\equiv \om_n\cdot r^n\quad \text{and}\quad \M(\D(S\on B_r(p)))\equiv n\cdot\om_n\cdot r^{n-1}\]
for all $r\leq r_0$.

We set $S_r:=S\on \bar B_r(p)$ and $\chi_r:=\spt(S_r)=\chi_S\cap\bar B_r(p)$.
 For $\la\in[0,1]$ denote by $c_\la:\bar B_{r_0}\to \bar B_{\la r_0}$
the map which contracts $p$ radial geodesics by the factor $\la$.
We claim that $(c_\la)_\# S_r$ is an area minimizer for every $r\leq r_0$. By slicing, 
\[\M(C_p(\D S_{r})\on(B_{r}(p)\setminus B_{\la r}(p)))\leq(1-\la^n)\cdot\frac{r}{n}\cdot\M(\D S_r)=(1-\la^n)\cdot\om_n\cdot r^n.\]
Since $c_\la$ is $\la$-Lipschitz, we have
\[\M((c_\la)_\# S_r)\leq\la^n \cdot\om_n\cdot r^n.\]
Because $C_p(\D S_{r})\on(B_{r}(p)\setminus B_{\la r}(p))+(c_\la)_\# S_r$ is a filling of $\D S_r$,
both inequalities above have to be equalities and $(c_\la)_\# S_r$ is an area minimizer with Euclidean mass growth relative $p$.

\bslem
For every $r\leq r_0$ and $\la\in[0,1]$ holds 
\[(c_\la)_\#S_r=C_p(\D(c_\la)_\# S_r).\]
\eslem

\proof
 Define
\[\varphi:\chi_r\times[0,1]\to X;\ (x,\la)\mapsto c_\la(x).\]
Since $\D C_p((c_\la)_\#S_r)=(c_\la)_\#S_r-C_p(\D(c_\la)_\# S_r)$, it is suffices to show $C_p((c_\la)_\#S_r)=0\in\bI_{n+1,c}(X)$.
Because
$\spt(C_p((c_\la)_\#S_r))\subset C_p(\chi_r)$,
it is enough to prove $\Ha^{n+1}(C_p(\chi_r))=0$.

The Lytchak-Rademacher theorem \cite[Theorem~1.6]{Ly_diff} implies the following.
Since $\chi_r$ is countably $\Ha^n$-rectifiable,
it has a linear tangent cone $T_{x}\chi_r\cong\R^n$ at almost all points. 
Moreover, since $\varphi$ is Lipschitz,
$\varphi$ is differentiable at almost all points $(x,\la)$
with a differential 
\[d\varphi(x,\la):T_{x}\chi_r\times\R\to T_{c_\la(x)}X\] 
which is linear onto its image. 
Suppose for contradiction that\\ 
$\Ha^{n+1}(C_p(\chi_r))>0$. Then the general area formula \cite[Theorem~8.2]{AK_rec} implies that there exists 
$\la_0\in(0,1)$ such that $d\varphi(x,\la_0)$ has rank $(n+1)$ for $x$ in a set of positive $\Ha^n$-measure.
Set $g:=(d_p)|_{\chi_S}$. Since $(c_{\la_0})_\#(S_r)$ is an area minimizer with Euclidean mass growth relative $p$,
the general coarea formula \cite[Theorem~9.4]{AK_rec} implies that
the coarea factor $C_1(d^S g_x)$ has to be $1$ almost everywhere on $\chi_r$, cf.~\cite[Section~9]{AK_rec}.
Let $x\in\chi_r$ be such that $\varphi$ is differentiable at $(x,\la_0)$ with linear differential of rank $(n+1)$
and $g$ has a tangential differential $d^S_x$ at $x$ \cite[Theorem~8.1]{AK_rec}. Then 
\[d\varphi(x,\la_0)(T_{x}\chi_r\times\R)\cong\R^{n+1}\subset T_{c_{\la_0}(x)}X.\]
Note that $\la\mapsto c_\la(x)$ is the geodesic from $x$ to $p$. Thus the image of $d\varphi(x,\la_0)$ contains $v_0:=\log_{c_{\la_0}(x)}(p)$.  
Let $\{v_1,\ldots,v_n\}$ be an orthonormal basis 
for $T_x \chi_r\cong\R ^n$ such that $v_1$ is the nearest point to $v_0\in T_x X$. Then, by the first variation formula, 
\[d^S g_x=(-\cos(\al_1),-\cos(\al_2),\ldots,-\cos(\al_n))=(-\cos(\al_1),0,\ldots,0)\] 
where $\al_i=\angle_{c_{\la_0}(x)}(v_0,v_i)$, $i\geq 1$. Here we have used that the directions $v_i$, $0\leq i\leq n$, lie in the image of 
$d\varphi(x,\la_0)$ which is an $(n+1)$-flat. Thus, $1=C_1(d^S g_x)=|\cos(\al_1)|$. Hence,
the tangent space $T_x\chi_r$ contains the direction $v_0\in T_x X$ and therefore  $d\varphi(x,\la_0)$
has rank at most $n$. Contradiction.
\qed

\bslem
For every $r\leq r_0$ and $\la\in[0,1]$ holds 
\[\D S_{\la r}=\D(c_\la)_\#S_r.\]
\eslem

\proof
Denote by $\pi_r:X\setminus B_r(p)\to S_r(p)$ the nearest point projection.
Note that $\D(S_r\on(X\setminus B_{\la r}(p)))=\D S_r-\D S_{\la r}$.
We have $(\pi_{\la r})_\#(S_r\on(X\setminus B_{\la r}(p)))=0\in\bI_{n,c}(X)$ since
\[\spt(\pi_{\la r})_\#(S_r\on(X\setminus B_{\la r}(p)))\subset\pi_{\la r}(\spt(\D S_r))\] 
and the latter is $\Ha^n$-negligible.
Since $c_\la=\pi_{\la r}$ on $\spt(\D S_r)$, we conclude
\[\D S_{\la r}=(\pi_{\la r})_\#(\D S_r)=(c_{\la})_\#(\D S_r)=\D(c_{\la})_\# S_r.\]
\qed

From the two sublemmas we conclude
\[S_{\la r}=C_p(\D S_{\la r})=C_p(\D(c_\la)_\# S_r)=(c_\la)_\# S_r.\]
By monotonicity and the mass control, the multiplicity function $\theta$ has to be equal to $1$ in a neighborhood of $p$ in $\chi_{r_0}$,
say on $\chi_{\la r_0}$ for appropriate $\la>0$.
We now show that $\theta\equiv 1$ on all of $\chi_{r_0}$.
Since $\spt(S_{\la r})\subset c_\la(\spt S_r)$, we have for any $r\leq r_0$,
\[
\om_n\cdot(\la r)^n=\M(S_{\la r})=\Ha^n(\chi_{\la r})
\leq\la^n\Ha^n(\chi_{r})\leq\la^n\M(S_r)=\om_n\cdot(\la r)^n
\]
Thus, $\|S\on B_{r_0}(p)\|=\Ha^n\on\chi_{r_0}$.
In particular, for all $r\leq r_0$ holds 
\[\Ha^n(\chi_S\cap B_r(p))=\om_n\cdot r^n\quad \text{and}\quad \Ha^{n-1}(\chi_S\cap S_r(p))= n\cdot\om_n\cdot r^{n-1}.\]
Since $\chi_S\cap \bar B_r(p)=\spt(S_r)\subset C_p(\spt\D S_r)=C_p(\chi_S\cap S_r(p))$ we conclude that the support of $S_r$ is conical
with respect to $p$, 
$\chi_S\cap \bar B_r(p)=C_p(\chi_S\cap S_r(p))$.
Moreover, $c_\la$ scales $\Ha^n$-measure on $\chi_S\cap B_{r_0}(p)$ exactly by the factor $\la^n$.
We claim that for every $x\in\spt(S)\setminus\spt(\D S)$ with $a=|p,x|$ and $r\leq r_0-a$ holds $\|S\|(\bar B_{r}(x))=\om_n\cdot r^n$.
Indeed, using the lower density bound at $x$, monotonicity and the effect of $c_\la$ on $\Ha^n$, we obtain
\begin{align*}
\om_n\cdot r^n&\leq \|S\|(B_r(x))=\la^{-n}\cdot \|S\|(c_\la(B_r(x)))\\
&\leq\la^{-n}\cdot \|S\|(B_{\la r}(c_\la(x)))\\
&\leq \la^{-n}\cdot\left(\frac{\la r}{r_0-\la a}\right)^n\cdot\om_n\cdot r_0^n
\stackrel{\la\to 0}{\longrightarrow} \om_n\cdot r^n.
\end{align*}
We infer $\chi_S\cap\bar B_r(x)=C_x(\chi_S\cap\bar B_r(x))$. Therefore $\chi_S\cap\bar B_{r_0}(p)$ is conical with respect to all of its points.
In particular, it is a closed convex subset of $X$. From the mass bound, it follows that  $\chi_S\cap\bar B_{r_0}(p)$ is isometric to a Euclidean $r_0$-ball.
\qed
\medskip

From the rigidity in Proposition~\ref{prop_mon}, we obtain:

\bcor\label{cor_mon}
Let $X$ be a locally compact CAT(0) space and $S \in \bI_{n,loc}(X)$ an area minimizer.
If $\Gi(S)=\om_n$, then $\spt(S)$ is an $n$-flat in $X$ and $S$ has constant multiplicity 1.
\ecor

\subsection{The asymptotic Plateau problem relative to a 2-flat}\label{subsec_relative}

We fix the following setup for this section. Let $X$ be a locally compact CAT(0) space with $1$-dimensional Tits boundary which contains a round $1$-sphere 
$\si\subset\tits X$ in its Tits boundary. We fix a 2-flat $F\subset X$ with $\tits F=\si$ and a base point $o\in F$. 
Moreover, $\al^+\subset\tits X$ denotes an embedded local geodesic with $\al^+\cap\si=\D\al^+$.
Our goal is to solve an asymptotic Plateau problem relative $F$. Informally speaking, we search for a minimal surface $S\subset X$
whose boundary lies in $F$ and which is asymptotic to $\al^+$. If the length of $\al^+$ is close to $\pi$, then we expect $S$
to behave almost like a flat half-plane orthogonal to $F$.
In order to achieve this, we employ doubling to reduce the problem to the ordinary asymptotic Plateau problem treated in \cite{KL_higher}. 
So we consider the double of $X$ along $F$:
\[\hat X=X^-\cup_F X^+\]
where $X^\pm$ are two isometric copies of $X$.  Note that $\tits \hat X=\tits X^-\cup_\si \tits X^+$.
In particular, $\dim(\tits \hat X)=1$. Let us denote by $\iota$ the natural isometry of $\hat X$ which interchanges $X^-$
and $X^+$. For subsets $A^+\subset X^+$ we set $A^-:=\iota (A^+)\subset X^-$. We call a current $S\in\bI_{2,loc}(\hat X)$
{\em symmetric}, if 
\[S=S^+-\iota_\#(S^+)\]
where $S^+:=S\on X^+$. Finally, we also double our boundary data, and define $\al:=\al^+\cup\iota(\al^+)\subset\tits X$.

\bdfn
For $S\in \bI_{n,loc}(\hat X)$  set $S^+:=S\on X^+$. Then $S^+$ is a {\em relative minimizer} if 
\[\M(S^+\on B^+)\leq \M(W^+)\]
whenever $B^+\subset X^+$ is a Borel set such that $S\on (B^+\cup B^-)\in\bI_{n,c}(\hat X)$ and $W^+\in\bI_{n,c}(X^+)$
satisfies $\spt(\D(S^+\on B^+-W^+))\subset F$.
\edfn

\blem\label{lem_flat_chain}
Let $S^+\in\cD_{2,loc}(X^+)$ be a current of locally finite mass and with $\spt(\D S^+)\subset F$.
Suppose that  $\|S^+\|(F)=0$ and that there is a sequence $t_k\to 0$, such that $S^+\on(X\setminus N_{t_k}(F))\in\bI_{2,loc}(X^+)$.
Then $\D S^+$ is a local flat chain in $F$, $\D S^+\in\cF_{1,loc}(F)$. Moreover, the double $S:=S^+-\iota_\#(S^+)$ is a locally integral cycle,
$S\in\bZ_{2,loc}(\hat X)$.
\elem

\proof
For a  point $p\in X^+$ we denote by $d_p$ the distance function to $p$.
To see that $\D S^+$ is a local flat chain in $F$, it is enough to prove the following claim.
For any $p\in\spt(\D S^+)$ there exists $r_0>0$ such that $\<S^+\on(X\setminus N_{t_k}(F)),d_p,r_0\>\in\bI_{1,c}(X^+)$ for every 
$k\in\N$, and $S^+\on\bar B_{r_0}(p)$ is a flat limit of integral chains $T_k\in\bI_{2,c}(X^+)$.

Indeed, if the claim holds, then, since the nearest point projection $\pi_F:X^+\to F$ is Lipschitz, $(\pi_F)_\#(S^+\on\bar B_{r_0}(p))\in\cF_{2,c}(F)$
and therefore $\D(\pi_F)_\#(S^+\on\bar B_{r_0}(p))=(\pi_F)_\#\D(S^+\on\bar B_{r_0}(p))\in\cF_{1,c}(F)$.
Then, either $p\notin\spt(\D S^+ -(\pi_F)_\#\D(S^+\on\bar B_{r_0}(p)))$ or, there exists a large $k\in\N$ such that 
$p\notin\spt(\D S^+ -(\pi_F)_\#\<S^+\on(X\setminus N_{t_k}(F)),d_p,r_0\>$. In any case, we see that $\D S^+$ is locally flat near $p$.

To prove the claim, we set $T_k:=S^+\on(X\setminus N_{t_k}(F))\in\bI_{2,loc}(X^+)$. By continuity of $S^+$ as an extended functional 
\cite[Theorem~4.4]{La_currents}, we have $T_k\to S^+\on(X^+\setminus F)$ weakly.
Since $\|S^+\|(F)=0$, we have $S^+\on(X^+\setminus F)=S^+$ \cite[Lemma~4.7]{La_currents}.
For $r>0$ we infer from \cite[Theorem~6.2 (3)]{La_currents}
\[\int_0^\infty \|\<T_k,d_p,s\>\|(B_r(p))\ ds\leq\|T_k\|(B_r(p))\leq\|S^+\|(B_r(p)).\]
By monotone convergence, we conclude
\[\int_0^\infty\lim\limits_{k\to\infty} \|\<T_k,d_p,s\>\|(B_r(p))\ ds\leq\|S^+\|(B_r(p)).\]
Now choose $r_0\in(0,r)$ such that $\<T_k,d_p,r_0\>\in\bI_{1,c}(X^+)$ for every $k\in\N$ and $\lim\limits_{k\to\infty} \|\<T_k,d_p,r_0\>\|<\infty$.
Set $T_k':=T_k\on\bar B_{r_0}(p)$. Because $\|S^+\|(F)=0$, we see $\lim\limits_{t\to 0}\|S^+\|(\bar B_{r_0}(p)\cap N_t(F))=0$. Then
for $k<l$ we have 
\[\cF(T'_k,T'_l)\leq (\|S^+\|+\lim\limits_{k\to\infty} \|\<T_k,d_p,r_0\>\|)(\bar B_{r_0}(p)\cap \bar N_{t_k}(F))\stackrel{k,l\to \infty}{\rightarrow}0.\]
Hence $(T'_k)$ is a Cauchy sequence with respect to flat distance.
Because flat convergence implies weak convergence, the limit coincides with $S^+\on\bar B_{r_0}(p)$. This confirms the claim and therefore $\D S^+\in\cF_{1,loc}(F)$.

By the deformation theorem \cite{FF_currents,W_defo}, we can now find a locally polyhedral chain $P\in\cP_{1,loc}(F)$ and a locally flat 
chain $V\in\cF_{2,loc}(F)$ such that $\D V=P-\D S^+$. Recall that since $V$ is top-dimensional in $F$, it is canonically identified with a locally integrable
function with integer values. In particular, $S^+ +V$ is locally integral. Hence $S=(S^+ +V)-\iota_\#(S^+ +V)$ is a locally integral cycle as claimed.
\qed

\bcor\label{cor_flat_chain}
Let $S\in\bZ_{2,loc}(\hat X)$ be symmetric, $S=S\on X^+-\iota_\#(S\on X^+)$.
Then the boundary of the restriction $S^+:=S\on X^+$ is a local flat chain in $F$, $\D S^+\in\cF_{1,loc}(F)$.
\ecor

\proof
Note that since $S$ is symmetric, we must have $\|S\|(F)=0$. Hence the claim follows from slicing and Lemma~\ref{lem_flat_chain}.
\qed

\blem\label{lem_minimizer}
There exists a minimizing $S\in\bZ_{2,loc}(\hat X)$  which is symmetric, and such that 
\[\Fi(S-C_o\al)=0\quad \text{ and }\quad \Gi(S)=\frac{\length(\al)}{2}.\] 
Moreover, the restriction $S^+:=S\on X^+$
is a relative minimizer.
\elem

\proof
The existence of a minimizer follows from \cite[Theorem~5.6]{KL_higher}. 
We only need to argue that we can choose it to be symmetric. As in \cite[Theorem~5.6]{KL_higher}, we choose a sequence of radii $r_k\to \infty$ and solve the Plateau Problem
for $C_o(\al)\on B_{r_k}(o)$ to obtain minimizers $S_k\in\bI_{2,c}(\hat X)$ with $\D S_k=\D (C_o(\al)\on B_{r_k}(o))$.
Note that since $C_o(\al)\on B_{r_k}(o)$ is symmetric, we may assume that $S_k$ is symmetric as well. Indeed, by Lemma~\ref{lem_flat_chain},  
$\tilde S_k:=(S_k\on X^+)-\iota_\#(S_k\on X^+)$ lies in $\bI_{2,c}(\hat X)$ and $\M(\tilde S_k)\leq 2\cdot\M(S_k\on X^+)\leq\M(S_k)$. 
Now the proof in  \cite[Theorem~5.6]{KL_higher} applies to produce a symmetric minimizer $S\in\bZ_{2,loc}(\hat X)$.
By \cite[Theorem~8.3]{KL_higher}, $S$ is $\F$-asymptotic to $C_o\al$ and has the required volume growth. That $S^+$ is a 
relative minimizer follows from symmetry of $S$.
\qed

\medskip

In Section~\ref{subsec_no_acc} we used flat half-planes orthogonal to $F$ to rule out accumulation of branch points at infinity.
Just like the relative minimizer $S^+$ serves as a substitute for a flat half-plane, the following lower bound on filling density serves as a substitute
for orthogonality to $F$.

From now on let $S\in\bZ_{2,loc}(\hat X)$ be a symmetric minimizer as in Lemma~\ref{lem_minimizer}.

\blem\label{lem_lower_dist_bound}
For every $p\in\spt(S)\cap F$ and $r>0$ holds
\[\F_{p,r}(S-F)\geq\frac{\pi}{3}.\]
\elem

\proof
Let $V\in\bI_{3,c}(\hat X)$ be such that $\spt(S-F-\D V)\cap B_r(p)=\emptyset$.
For almost all $s\in[0,r]$ we have $\<S,d_p,s\>-\<F,d_p,s\>=\D\<V,d_p,s\>$ and all occurring slices are integral. Since $S$ is a symmetric minimizer, we must have
\[\|\<V^\pm,d_p,s\>\|(B_s(p))\geq\|S^\pm(B_s(p))\|.\]
We obtain $\|S^\pm(B_s(p))\|\geq\frac{\pi}{2}\cdot s^2$ from Proposition~\ref{prop_mon}. Since 
\[\M(V)\geq 2\cdot\M(V^+)\geq 2\cdot\int_0^r \M(\<V^+,d_p,s\>)\ ds,\]
the claim follows by integration.
\qed

\blem\label{lem_rel_visibility}
For every $\eps>0$ there exists $r_\eps>0$
such that for all $r\geq r_\eps$ holds
\[\spt(\D S^+)\cap B_r(o)\subset N_{\eps r}(C_o(\D\al^+)).\]
\elem

\proof
Suppose for contradiction that there exists $\eps_0>0$ and a sequence $x_k\in\spt(\D S^+)$
with $|o,x_k|\to\infty$ but $|x_k,C_o(\D\al^+)|\geq\eps_0\cdot|o,x_k|$. 
Since $\spt(\D S^+)\subset\spt(S^+)\subset\spt(S)$, we infer from \cite[Theorem~8.1]{KL_higher}, that there exist points $y_k\in C_o(\al)$ with $\frac{|x_k,y_k|}{|x_k,o|}\to 0$.
After passing to a subsequence, we find an element $\eta\in\al$ such that $\frac{|y_k,o\eta|}{|y_k,o|}\to 0$.
Hence $\frac{|x_k,o\eta|}{|x_k,o|}\to 0$. But since $x_k\in F$, we have $\eta\in\D\al^+$. Contradiction.
\qed

\blem\label{lem_intersection}
Let $\zeta$ and $\hat\zeta$ be antipodes on $\si$
which separate $\D\al^+$. Then $\spt(\D S^+)$ intersects every complete geodesic $g_0\subset F$ with $\geo g_0=\{\zeta,\hat\zeta\}$.
\elem

\begin{center}
\includegraphics[scale=0.5,trim={-3cm 0cm 0cm 0cm},clip]{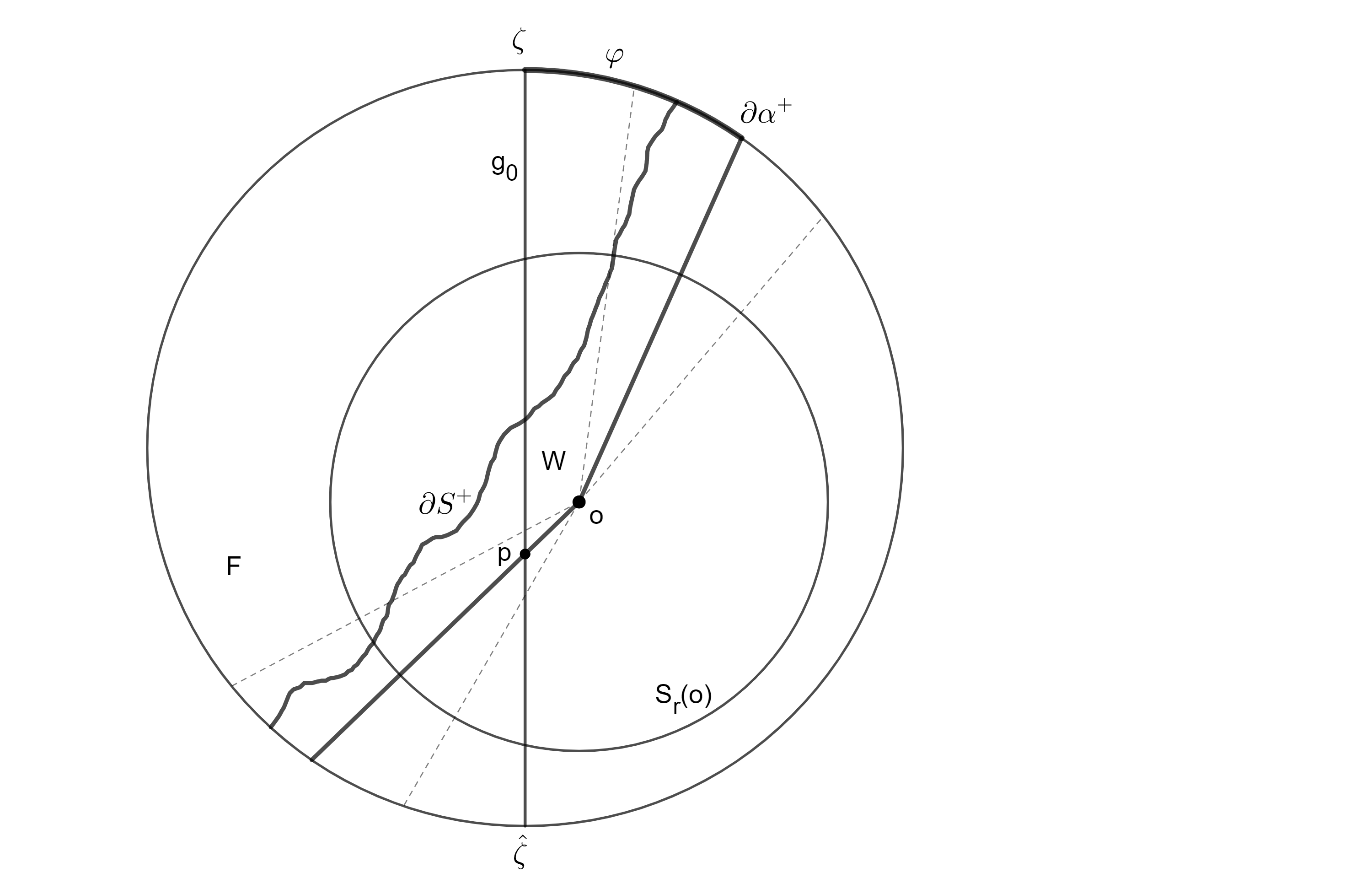}
\end{center}

\proof
By Lemma~\ref{lem_flat_chain}, $\D S^+$ is a local flat chain in $F$.
Note that $\Fi(S-C_o(\al))=0$ implies $\Fi(\D S^+-C_o(\D\al))<\infty$. 
Indeed, for $\eps_0>0$ choose $r>0$ and $V\in\bI_{3,c}(\hat X)$ such that 
$\spt(S-C_o\al-\D V)\cap B_r(o)=\emptyset$ and $\M(V)\leq\eps_0\cdot r^3$.
Denote by $d_F$ the distance function to $F$ in $X^+$ and by $\pi_F:\hat X\to F$ the nearest point projection.
By the coarea inequality, we find $s\in(0,\frac{r}{2})$ such that $\M(\<V,d_F,s\>)\leq\frac{2}{r}\cdot\M(V)\leq 2\eps_0\cdot r^2$
and $\spt(\<S-C_o\al,d_F,s\>-\D\<V,d_F,s\>)\cap B_r(o)=\emptyset$, and all occurring slices are integral.
Set $T:=(S^+-C_o\al^+)\on\bar N_s(F)-\D\<V,d_F,s\>\in\bI_{2,c}(X^+)$. Then 
$\spt(\D S^+-C_o\D\al^+-\D(\pi_F)_\# (T))\cap B_{\frac{r}{2}}(o)=\emptyset$ and $\M((\pi_F)_\# (T))\lesssim r^2$.

In particular, $\D S^+$ and $C_o(\D\al^+)$ are homologous as local flat chains.
Denote by $W\in\cF_{2,loc}(F)$ the canonical filling of $\D S^+ -C_o(\D\al^+)$. By Lemma~\ref{lem_rel_visibility},
the support of $W$ lies sublinearly close to $C_o(\D\al^+)$.
Let $\varphi$ be a positive lower bound for the distance between the sets $\{\zeta,\hat\zeta\}$ and $\D\al^+$ in $\tits F$.
Choose $\eps\leq\frac{\varphi}{2}$ and then choose $r_\eps>0$ according to Lemma~\ref{lem_rel_visibility}.
By assumption, $C_o(\D\al^+)$ intersects $g_0$ transverse in exactly one point, say $p\in F$.
Now we can find $r\geq \max\{r_\eps,2|o,p|\}$ such that all points in $\spt(\<W,d_o,r\>)$ have distance larger than $\eps\cdot r$ from the geodesic $g_0$.
 Since $\D(W\on B_r(o))=\<W,d_o,r\>+(\D S^+-C_o(\al^+))\on B_r(o)$ is a flat 1-cycle with compact support,
$\spt(\D S^+)$ has to intersect $g_0$ in $B_r(o)$. 
\qed

\subsection{Producing flat half-planes as relative minimizers}\label{subsec_rel_asym}

Now we have all the necessary ingredients to prove:

\begin{namedlemma}[Half-Plane]
Let $X$ be a locally compact CAT(0) space with $1$-dimensional Tits boundary.  Let $F\subset X$ be a 2-flat with $\geo F=\si$.
Further, let $g$ be an axis of an axial isometry $\ga$ and assume $\geo g\subset\si$.
Let $\xi^\pm\in\si$ be antipodes disjoint from $\geo g$.
Suppose that there is a sequence of local geodesics $\al^+_k\subset \tits X$ such that 
\begin{itemize}
	\item $\D\al^+_k=\{\xi^-_k,\xi^+_k\}$ and $\xi^\pm_k\to\xi^\pm$ with respect to the Tits metric;
	\item the length of $\al_k^+$  converges to $\pi$ as $k\to\infty$.
\end{itemize}
Then there exists a 2-flat $\tilde F\subset P(g)$, and a flat half-plane $\tilde H$ with boundary $\tilde h:=\D \tilde H\subset\tilde F$ such that the following properties hold.
\begin{enumerate}
	\item $g\subset N_a(\tilde F)$ if $g\subset N_a(F)$ for $a>0$;
	\item $|\geo \tilde h,\geo g|=|\{\xi^-,\xi^+\},\geo g|$;
	\item $\geo\tilde H\cap\geo\tilde F=\geo\tilde h$.
\end{enumerate}
\end{namedlemma}

\begin{center}
\includegraphics[scale=0.5,trim={-1cm 0cm 0cm 0cm},clip]{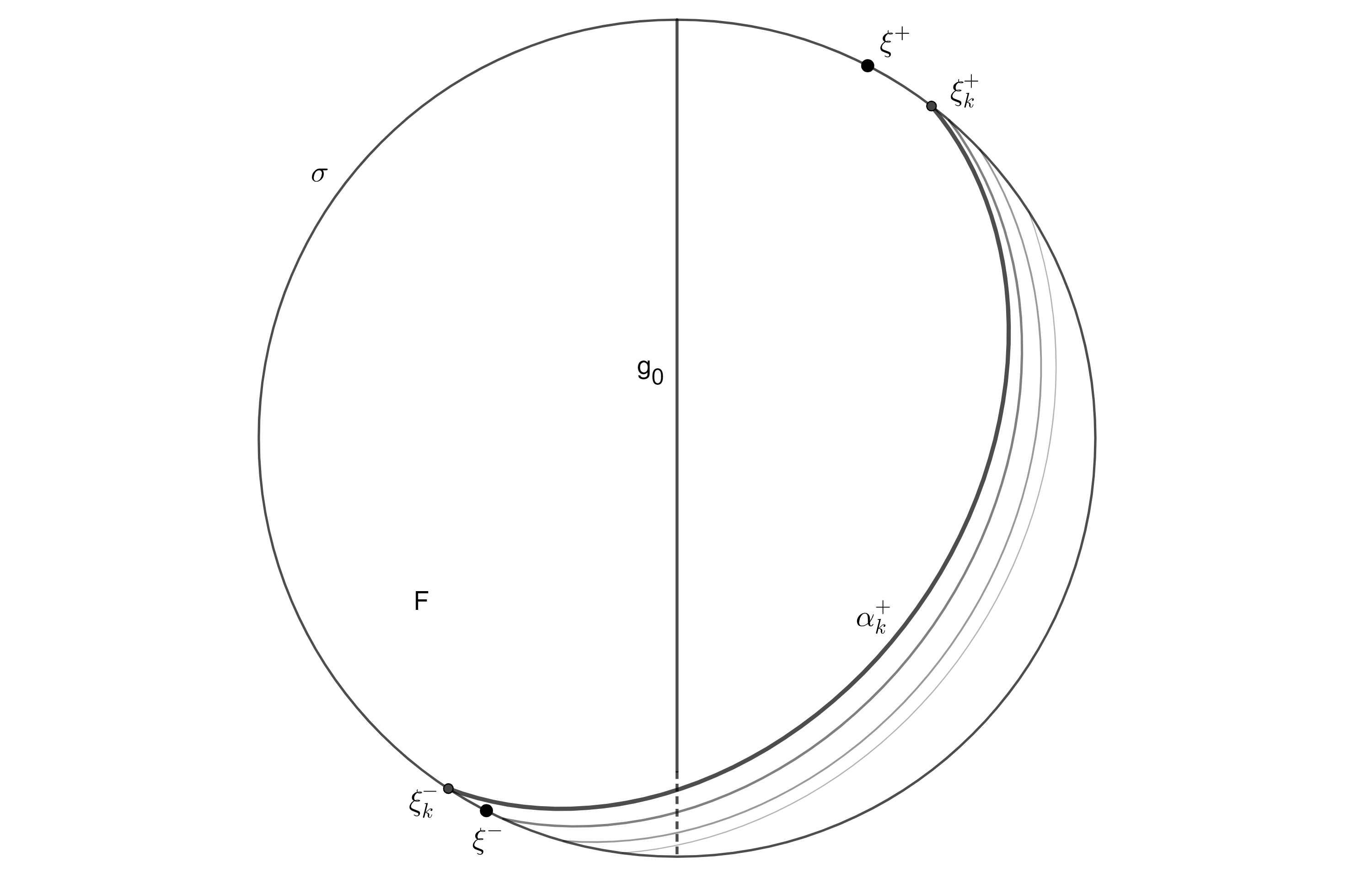}
\end{center}

\proof
Choose a  base point $o\in F$. 
As before, we consider the double of $X$ along $F$, $\hat X=X^-\cup_F X^+$
and put $\al_k:=\al_k^+\cup \iota(\al_k^+)$.
By Lemma~\ref{lem_minimizer}, we find a symmetric minimizer $S_k\in\bZ_{2,loc}(\hat X)$ with $\Fi(S_k-C_o(\al_k))=0$
and $\Gi(S_k)=\frac{\length(\al_k)}{2}$. Let $g_0$ be a complete geodesic in $F$ parallel to $g$. By Lemma~\ref{lem_intersection},
there exists  a point $p_k\in\spt(S_k)\cap g_0$. 

Now let us choose numbers $m_k\in\Z$ such that the points $\ga^{m_k}(p_k)$ stay in a fixed compact set. Put $\ga_k:=\ga^{m_k}$.
After passing to subsequences, we have convergences with respect to pointed Hausdorff topology 
\[\ga_k(C_{p_k}(\al^+_k),p_k)\to(\tilde H,p_\infty)\quad\text{ and }\quad\ga_k (F,p_k)\to(\tilde F,p_\infty)\]
where $\tilde H$ is a flat half-plane, $\tilde F$ is a 2-flat, and $\tilde h:=\D\tilde H\subset\tilde F$.
We will now show that $\tilde F$ and $\tilde H$ satisfy the three required properties.
\medskip

\noindent (1) Since $\ga$ preserves $g$, we have $g\subset N_a(\tilde F)$ for $a>0$ whenever $g\subset N_a(F)$.
\medskip
	
\noindent (2) Let $\geo\tilde h=\{\tilde\xi^-,\tilde\xi^+\}$. 
Since $\xi^\pm_k\to\xi^\pm$ with respect to the Tits metric, and since $\ga$ preserves $g$, we have
$|g(\pm\infty),\ga_k\xi_k^+|=|g(\pm\infty),\xi_k^+|$ and $|g(\pm\infty),\ga_k\xi_k^-|=|g(\pm\infty),\xi_k^-|$.
Hence the statement follows from lower semi-continuity of the Tits metric with respect to the cone topology
\cite[Lemma~2.3.1]{KleinerLeeb}.
\medskip

\noindent (3) Note that if this fails, then  we must have $\tilde H\subset \tilde F$ since $\D\tilde H\subset\tilde F$.

Denote by $\tilde F^\pm\subset\tilde F$ the two half-planes determined by $\D\tilde H$.
Further, denote by $\si_k^\pm\subset\si$ the two closed arcs determined by $\{\xi_k^-,\xi_k^+\}$.
Then, as currents, we can write $F=F_k^+-F_k^-$ where $F_k^\pm=C_{p_k}(\si_k^\pm)\in\bI_{2,loc}(\hat X)$.
From \cite[Proposition~4.5]{KL_higher}, we conclude that for every $\eps>0$
there exists $a_\eps'>0$ such that $\F_{p_k,r}(S_k-C_{p_k}(\al_k))\leq\frac{\eps}{3}$ for all $r\geq a_\eps'$.
On the other hand, Lemma~\ref{lem_lower_dist_bound} implies $\F_{p_k,r}(S_k-F)\geq\frac{\pi}{3}$ for all $r>0$.
Hence, $\F_{p_k,r}(F-C_{p_k}(\al_k))\geq\frac{\pi-\eps}{3}$ and thus, by symmetry of  $C_{p_k}(\al_k)$,  for $r\geq a_\eps'$, 
\[\F_{p_k,r}(F_k^\pm-C_{p_k}(\al^+_k))\geq\frac{\pi-\eps}{6}.\tag{$\star$}\] 
Since $\length(\al_k^+)\to\pi$, we have the uniform  mass bound 
\[(\|C_{p_k}(\al^+_k)\|+\|\D C_{p_k}(\al^+_k)\|)(B_r(p_k))\leq \pi r^2+2r.\]
Using \cite[Theorem~2.3]{KL_higher}, we may assume that $(\ga_k)_\#(C_{p_k}(\al^+_k))$ converges to a limit $C_\infty\in\bI_{2,loc}(X^+)$ with respect to local flat topology. By \cite[Proposition~2.2]{W_compact}, the support of $C_\infty$ is contained in $\tilde H$. 
Because $(\ga_k)_\#\D(C_{p_k}(\al^+_k))\to\D\tilde H$ weakly, we have $\D C_\infty=\D\tilde H$. But $C_\infty$
is a top-dimensional locally integral current in $\tilde H$, hence $C_\infty=\tilde H$ \cite[Theorem~7.2]{La_currents}.

Now suppose for contradiction, that $\tilde H=\tilde F^+$. Set $T_k=F_k^+-C_{p_k}(\al^+_k)\in\bZ_{2,loc}(X^+)$. Because $(\ga_k)_\#F_k^+\to\tilde F^+$ locally flat, we infer $(\ga_k)_\#(T_k)\to 0$ locally flat. Since $\D T_k=0$ for all $k\in\N$,
\cite[Proposition~3.2]{LaWe_loc} implies $\F_{p_\infty,r}((\ga_k)_\#(T_k))\to 0$ for every $r>0$.
Since $\ga_k(p_k)\to p_\infty$, this contradicts ($\star$). Hence $\tilde H$ cannot be contained in $\tilde F$ and the proof is complete.
\qed

\bibliographystyle{alpha}
\bibliography{rr_III}

\Addresses

\end{document}